\documentclass[11pt,fleqn]{article}
\usepackage{setspace}
\RequirePackage[OT1]{fontenc}
\RequirePackage{amsthm,amsmath}
\RequirePackage[round]{natbib}
\RequirePackage[colorlinks,citecolor=blue,urlcolor=blue]{hyperref}
\usepackage{pgf,tikz}\usepackage{mathrsfs}\usetikzlibrary{arrows}

\usepackage{float}
\usepackage{fullpage}
\usepackage{amsmath}
\usepackage{amssymb}
\usepackage{graphicx}
\usepackage{dsfont}
\makeatletter
\def\captionof#1#2{{\def\@captype{#1}#2}}
\makeatother
\def\1{\mbox{\bf 1}}
\def\R{\mathbb{R}}

\def\N{\mathbb{N}}
\def\P{\mathbb{P}}
\def\E{\mathbb{E}}
\def\L{\mathbb{L}}

\def\R{\mathbb{R}}

\def\Z{\mathbb{Z}}

\def\c{\mbox{Cov}}

\newtheorem{theo}{Theorem}

\newtheorem{prop}{Proposition}

\newtheorem{Def/Prop}{Definition-Proposition}

\newcounter{exos}
\renewcommand\theexos{\arabic{exos}}

\newcounter{prob}
\renewcommand\theprob{\arabic{prob}}

\newcommand{\footremember}[2]{%
    \footnote{#2}
    \newcounter{#1}
    \setcounter{#1}{\value{footnote}}%
}
\newcommand{\footrecall}[1]{%
    \footnotemark[\value{#1}]%
}

\begin{document}
\author{Zinsou Max Debaly \footremember{affil}{CREST-ENSAI, UMR CNRS 9194, Campus de Ker-Lann, rue Blaise Pascal, BP 37203, 35172 Bruz cedex, France.}\and
Lionel Truquet \footrecall{affil}\footnote{This work was funded by CY Initiative of Excellence
(grant "Investissements d'Avenir" ANR-16-IDEX-0008),
Project "EcoDep" PSI-AAP2020-0000000013.}}

\title{Iterations of dependent random maps and exogeneity in nonlinear dynamics}
\date{}
\maketitle

\begin{abstract}
\noindent
We discuss the existence and uniqueness of stationary and ergodic nonlinear autoregressive processes when exogenous regressors are incorporated into the dynamic.
To this end, we consider the convergence of the backward iterations of dependent random maps. In particular, we give a new result when the classical condition of contraction on average is replaced with a contraction in conditional expectation. Under some conditions, we also discuss the dependence properties of these processes using the functional dependence measure of Wu (2005) that delivers a central limit theorem giving a wide range of applications. Our results are illustrated with CHARN models, GARCH processes, count time series, binary choice models and categorical time series for which we provide many extensions of existing results.
\end{abstract}
\vspace*{1.0cm}

\footnoterule

\section{Introduction}
Among the various contributions devoted to time series analysis, theoretical results justifying stationarity and ergodicity properties of some standard stochastic processes when exogenous covariates are incorporated in the dynamic are rather scarce. A notable exception concerns linear models, such as VARMA processes, for which such properties are a consequence of the linearity. See for instance \citet{Luthp}, a standard reference for multivariate time series models.
Moreover, linear models represent a very simple setup for discussing various exogeneity notions found in the literature. See for instance \citet{Engle}.
For nonlinear dynamics, a few contributions consider the problem of exogenous regressors. 
For general GARCH type processes, \citet{Francq} recently studied stationarity conditions when the noise and the covariate process form a stationary process. \citet{Cav} considered a Poisson autoregressive process with exogenous regressors (PARX models), under a Markov chain assumption for the covariate process. 
\citet{deJ} consider the case of dynamic binary choice models and provide results about stationarity and mixing properties of a $0/1-$valued time series which is autoregressive and defined conditionally on some exogenous regressors.  
\citet{FT1} studied stationarity and ergodicity of general categorical time series defined conditionally on a strictly exogenous covariate process.

 In this paper, we give general results for getting stationarity, ergodicity and stochastic dependence properties
for general nonlinear dynamics defined in terms of iterations of random maps. For simplicity, we explain our setup with the following example which represents the basis for studying other processes. Let us consider the following model
\begin{equation}\label{eq::discuss}
X_t=F\left(X_{t-1},Z_{t-1},\varepsilon_t\right),\quad t\in\Z,
\end{equation}
where $(Z_t)_{t\in\Z}$ is a covariate process and $\left(\varepsilon_t\right)_{t\in\Z}$ 
a noise process. One can note that $X_t=f_t(X_{t-1})$ for the random function defined by 
$f_t(x)=F\left(x,Z_{t-1},\varepsilon_t\right)$. The sequence $(f_t)_{t\in\Z}$ is a sequence of dependent random maps even if the $\varepsilon_t'$s are i.i.d. because typically the $Z_t'$s exhibit temporal dependence.
A key point for getting existence of a stationary solution in (\ref{eq::discuss}) is to control the behavior of the backward iterations
$\left\{f_t\circ f_{t-1}\circ\cdots\circ f_{t-n}(x): n\geq 1\right\}$. 
The convergence of such iterations of random maps has been extensively studied in the independent case.
In this case, the process $(X_t)_{t\in\Z}$ is a Markov chain.
We defer the reader to \citet{Letac} and \citet{Diac} for seminal papers on iterated independent random maps and to \citet{Wu2} for additional results useful in a time series context. The last contribution is particularly interesting for getting existence of some moments for the marginal $X_t$ and also some dependence properties for the process $(X_t)_{t\in\Z}$ that are often needed for statistical applications. All these contributions use average contraction conditions and the interested reader is referred to the interesting survey of \citet{Stenflo} for an overview of the available results. There also exist some contributions studying the more general case of iterated stationary random maps $(f_t)_{t\in\Z}$. For instance,
\citet{Borovkov} gives many results for studying what he calls stochastically recursive sequences, when the independence assumption is removed. See also \citet{Iosif} for a survey of some available results.
The results obtained in the dependent case are based on Lyapunov type exponents and the convergence of the backward iterations is only studied almost surely. 
We recall the following result which can be found in \citet{Elton}  (see also \citet{Iosif}, Theorem $6.2$) and which generalizes 
a widely known result given in \citet{Brandt} or \citet{BP} for iterations of affine random maps.

To do so, we introduce some notations and conditions. We assume that $f_t:E\rightarrow E$ are random Lipschitz functions where $E$ denotes a locally compact Polish space endowed with a metric $d$. We define the Lipschitz constant of a measurable function $g:E\rightarrow E$ by
$$c(g):=\sup_{x\neq y\in E}\frac{d\left(g(x),g(y)\right)}{d(x,y)}.$$
Moreover, for any integers $s<t$, we set $f_s^t=f_t\circ\cdots\circ f_s$. For a positive real number $x$, we set $\log^{+}(x)=\log(x)$ if $x\geq 1$ and $0$ otherwise.
\begin{theo}\label{bibli}
Assume that the process $\left((Z_t,\varepsilon_t)\right)_{t\in\Z}$ in (\ref{eq::discuss}) is stationary and ergodic. Assume further that $\E\left[\log^{+}c(f_0)\right]<\infty$ and $\E\left[\log^{+}d(x_0,f_0(x_0))\right]<\infty$ for some point $x_0\in E$. 
\begin{enumerate}
\item
There exists a constant $\chi\in \R\cup\{-\infty\}$ called Lyapunov exponent and such that 
$$\lim_{n\rightarrow \infty}\frac{1}{n}\log c\left(f_1^n\right)=\chi\mbox{ a.s}.$$
Moreover 
$$\chi=\lim_{n\rightarrow \infty}\frac{1}{n}\E\left[\log c\left(f_1^n\right)\right]=\inf_{n\geq 1}\frac{1}{n}\E\left[\log c\left(f_1^n\right)\right].$$
\item
If the constant $\chi$ is negative, then the almost sure limit $f_{-\infty}^t=\lim_{k\rightarrow \infty}f_{t-k}^t(x)$ exists for any $x\in E$ and does not depend on $x$. Setting $X_t=f_{-\infty}^t$, the process $(X_t)_{t\in\Z}$ is stationary and ergodic and satisfies the recursions (\ref{eq::discuss}). Moreover, $(X_t)_{t\in \Z}$ is the unique stationary process satisfying (\ref{eq::discuss}).
\end{enumerate}
\end{theo}
The affine random maps version of this result has been applied recently by \citet{Francq} for studying stationarity of asymmetric power GARCH processes. For nonlinear random maps, Theorem \ref{bibli}
is less known in the time series literature. In this paper, we will make use of Theorem \ref{bibli} for defining a general class of categorical time series with exogenous covariates. In particular, we will see in Section \ref{binary} how Theorem \ref{bibli} can be applied to binary time series and lead to an improvement of a result of \citet{deJ}. 

However, the result presented above has several limitations. 

\paragraph{1.} First, it requires the random maps $f_t$ to be almost surely Lipschitz. Such a property is not always valid, for instance for the Poissonian autoregressions discussed in Section  \ref{Poissss}. When there are no exogenous covariates, \citet{davis} studied integer-valued time series by using a different contraction result, developed by \citet{Wu2}.  
 
\paragraph{2.} Existence of some moments for the marginal distributions that are sometimes necessary for statistical applications cannot be obtained directly from this result.
 
\paragraph{3.} For autoregressions with several lags, it is not straightforward to get an explicit condition on the parameters of the model to ensure that $\chi<0$.

To overcome these drawbacks, we will adapt the approach used by \citet{Wu2} for independent random maps to the case of dependent random maps. 
Our main result, see Theorem \ref{main} and its extension Theorem \ref{autoreg}, is obtained by replacing the usual contraction on average condition by a contraction in conditional expectation.
The assumptions that we use are very simple to check and the proof of our main result is straightforward but its merit is to provide an elegant way for presenting a general approach which encompasses most of the previous attempts to include exogenous regressors in nonlinear dynamics. For strictly exogenous regressors, i.e. the processes $(Z_t)_{t\in \Z}$ and $\left(\varepsilon_t\right)_{t\in\Z}$ are independent, we also provide an additional result, see Theorem \ref{relax}, with weaker assumptions.
In the context of Theorem \ref{main} and Theorem \ref{autoreg}, we will then discuss how to control the functional dependence measure of $(X_t)_{t\in\Z}$ introduced by \citet{Wu}, a dependence notion which is an alternative to the standard strong mixing condition and which can be more easily checked for iterations of contracting random maps.   
Let us mention that even in the independent case, mixing properties of the process $(X_t)_{t\in\Z}$ require restrictive assumptions on the noise distribution, otherwise such properties may fail. We refer the reader to the standard textbook of \citet{Doukhan(1994)}, section $2.4$ for mixing properties of iterations of independent random functions and to \citet{Andrews} for a famous counterexample of a non strongly mixing sequence defined via iterations of random maps.
In the dependent case, as in (\ref{eq::discuss}), getting usual strong mixing properties 
seems to be harder because the process $(X_t)_{t\in\Z}$ does not have a Markov structure in general and the criteria for getting mixing properties of Markov chains are useless.   

This paper is mainly motivated by dynamics of type (\ref{eq::discuss}) with covariates that are not necessarily strictly exogenous, assuming that at any time $t$,
the noise $\varepsilon_t$ is independent from the past information $\sigma\left((Z_s,\varepsilon_s):s\leq t-1\right)$.
The term predetermindness is sometimes used in the literature.
This independence assumption is substantially weaker than the independence between the two processes $\left(\varepsilon_t\right)_{t\in\Z}$ and $\left(Z_t\right)_{t\in\Z}$. The latter independence condition implies strict exogeneity, a notion initially defined by \citet{Sims} and extended to general models by \citet{Chamb}. Strict exogeneity is useful for deriving the conditional likelihood of the $X_t'$s conditionally on the $Z_t'$s. 
However, strict exogeneity is a rather strong assumption. Under additional regularity conditions on the model, \citet{Chamb} has shown that this assumption is equivalent to the non Granger-causality, i.e.
$Z_t$ is independent of $(X_s)_{s\leq t}$ conditionally on $(Z_s)_{s\leq t-1}$. It roughly means that the covariate process $(Z_t)_{t\in\Z}$ evolves in a totally autonomous way. 
In contrast, our exogeneity condition allows general covariates of the form $Z_t=H\left(\eta_t,\eta_{t-1},\ldots\right)$ with $H$ a measurable function and a sequence $\left((\eta_t,\varepsilon_t)\right)_{t\in\Z}$ of i.i.d. random vectors, $\varepsilon_t$ being possibly correlated with $\eta_t$. The error $\varepsilon_t$ 
can then still have an influence on future values of the covariates.
For linear models, the two technical independence conditions discussed above between the noise and the covariate processes are often used as a distinction between weak and strict exogeneity. See for instance \citet{Luthp}, Section $10.2$.
Let us mention that there exist additional concepts of exogeneity that are introduced and discussed in \citet{Engle}, in particular a notion of weak exogeneity.  However, this notion is related to the estimation of a specific parameter of the conditional distribution for the bivariate process $(X_t,Z_t)$ and it is necessary to specify the joint dynamic of the process. Since we do not want to consider specific dynamics for the covariate process, we will not use it in this paper.  
Inclusion of exogenous regressors motivates our approach which is based on conditional average contraction conditions. But our results can be also applied without referring to these concepts of exogeneity, i.e. when $\left((Z_{t-1},\varepsilon_t)\right)_{t\in\Z}$ is a general stationary and ergodic process in (\ref{eq::discuss}). However, in the latter case, a closed form expression for the conditional distribution of $X_t$ given $\mathcal{F}_{t-1}$ cannot be obtained directly from the recursions (\ref{eq::discuss}).  Our contribution is then the first one presenting a general framework for inclusion of covariates in nonlinear dynamics.  Our results can be applied to any statistical procedure which require either ergodic properties or the use of some limit theorems developed from the notion of functional dependence introduced by \citet{Wu}.  
Since the existing literature already contains many asymptotic results of this type, we do not discuss specific applications.
  
The rest of this paper is organized as follows. In Section \ref{general}, we give our main results for defining stationary and ergodic solutions for recursions of type (\ref{eq::discuss}). In Section \ref{mix}, we study weak dependence properties of the process using the functional dependence measure of \citet{Wu}. Many examples of nonlinear time series models satisfying our assumptions are given in Section \ref{examples} and we revisit some nonlinear dynamics discussed recently in the literature but we also consider new ones. 
A conclusion is given in Section \ref{Conc}. The proofs of our results are postponed to the last section of the paper.

\section{General results}\label{general}
In this section, we state several results for controlling the convergence of the backward iterations in some $\L^p$ spaces.
We recall that for a random variable $X$ and a real number $p\geq 1$, the quantity $\Vert X\Vert_p=\E^{1/p}\left(\vert X\vert^p\right)$
is called the $\L^p-$norm of $X$.  
Now let $(f_t)_{t\in\Z}$ be a sequence of random maps defined on a Polish space $(E,d)$ and taking values in the same space. We assume for convenience that $f_t=F\left(\cdot,\zeta_t\right)$ where 
$(\zeta_t)_{t\in\Z}$ is a stochastic process taking values in another Polish space $E'$ and $F:E\times E'\rightarrow E$ is a measurable map. In connection with our initial example (\ref{eq::discuss}), we have $\zeta_t=\left(Z_{t-1},\varepsilon_t\right)$. 
In the latter case, we will assume throughout the paper that $E'=E'_1\times E'_2$ where $E_1'$ is a Borel subset of $\R^e$ and $E_2'$ is another Polish space.

For $s<t$, we set $f_s^t=f_t\circ f_s^{t-1}$ with the convention $f_t^t=f_t$ and $f_t^{t-1}(x)=x$.
Moreover, we consider a filtration $\left(\mathcal{F}_t\right)_{t\in\Z}$ for which $\left(\zeta_t\right)_{t\in\Z}$ is adapted.

\subsection{Conditional contraction on average}
We first give a general and useful result for getting a convergence in some $\L^p$ spaces.  
For some real numbers $p\geq 1$, $L > 0$, $\kappa \in (0,1)$ and an integer $m  \geq 1$, we consider the two following assumptions.

\paragraph{A1}
There exists $x_0\in E$ such that $\sup_{t\in\Z}\E\left[d^p\left(f_t(x_0),x_0\right)\right]<\infty$.

\paragraph{A2}
For every $t\in\Z$, almost surely, the following inequalities hold for every $(x,y)\in E^2$.   
$$\E\left[d^p\left(f_t(x),f_t(y)\right)\vert \mathcal{F}_{t-1}\right]\leq L^pd^p(x,y)\mbox{ and }\E\left[d^p\left(f^{t+m-1}_t(x),f^{t+m-1}_t(y)\right)\vert \mathcal{F}_{t-1}\right]\leq \kappa^pd^p(x,y).$$

\begin{theo}\label{main}
Suppose that Assumptions {\bf A1-A2} hold. 
\begin{enumerate}
\item
For every $(x,t)\in E\times \Z$, there exists a $E-$valued random variable $X_t(x)$ such that 
$$\sup_{t\in\Z}\Vert d\left(X_t(x),x_0\right)\Vert_p<\infty,\quad \sup_{t\in\Z}\Vert d\left(f_{t-s}^t(x),X_t(x)\right)\Vert_p=O\left(\kappa^{s/m}\right).$$
Moreover the sequence $\left(f^t_{t-s}(x)\right)_{s\geq 0}$ converges almost surely to $X_t(x)$.
\item
For $x\neq y$, we have $\P\left(X_t(x)\neq X_t(y)\right)=0$. We then set $X_t=X_t(x)$. 
\item
The process $\left((X_t,\zeta_t)\right)_{t\in\Z}$ is stationary and also ergodic if the process $(\zeta_t)_{t\in\Z}$ is itself stationary and ergodic.
\item
If $(Y_t)_{t\in\Z}$ is a non-anticipative process (i.e. $Y_t\in \mathcal{F}_t$) such that $Y_t=f_t\left(Y_{t-1}\right)$ for $t\in\Z$ and 
$\sup_{j\in\Z}\E\left[d^p\left(Y_j,x_0\right)\right]<\infty$, then $Y_t=X_t$ a.s.
\end{enumerate}
\end{theo}

\paragraph{Notes} 
\begin{enumerate}
\item
The bounds given in Assumption {\bf A2} are required to hold for all $(x,y)\in E^2$ at the same time.
Since a conditional expectation is only unique up to a set with measure $0$ for the probability measure $\P$, the bound given in Assumption {\bf A2} has to be understood in term of regular conditional distribution, i.e. there exists a regular version of the conditional distribution of $\left(\zeta_t,\ldots,\zeta_{t+m-1}\right)$ given $\mathcal{F}_{t-1}$. On Polish spaces, a regular version always exists. See \citet{Kall}, Chapter $5$.   
\item
When $\zeta_t=\left(Z_{t-1},\varepsilon_t\right)$ forms a stationary process and $\mathcal{F}_t=\sigma\left((Z_j,\varepsilon_j): j\leq t\right)$, Theorem \ref{main} guarantees existence and uniqueness of a stationary process possessing a moment of order $p$ and solution of (\ref{eq::discuss}).
However, stationarity of the covariate/error process is not required for applying this result. In particular, when $(\zeta_t)_{t\in\Z}$ is non stationary, one can still define solutions of (\ref{eq::discuss}) provided that {\bf A1-A2} are satisfied. In general, these solutions will be non stationary and the classical law of large numbers is not valid. In this case, studying the asymptotic properties of some classical inferential procedures such as conditional likelihood estimation requires a specific analysis. 
 
\item 
Setting $p=1$ and 
$d(x,y)=\vert x-y\vert^o$ for $x,y\in E=\R$ and some $o\in (0,1)$, one can consider stochastic recursions (\ref{eq::discuss}) with  heavy-tailed covariate processes.
 
\item
When the process $(\zeta_t)_{t\in\Z}$ is stationary and the recursions are initialized at time $t=0$ with a given state $x\in E$,
the probability distribution of the forward iterations $f_1^t(x)$ coincides with the probability distribution of the backward iterations $f_{-t+1}^0(x)$. Since $\lim_{t\rightarrow \infty} f_{-t+1}^0(x)=X_0$ a.s.,  $f_1^t(x)$ converges in distribution to $X_0$. The same property holds true when the iterations are initialized with a random variable $\overline{X}_0$ independent from $\left(\zeta_t\right)_{t\geq 1}$. The main interest of the convergence of the backward iterations is to define the good random initialization $X_0=f_{-\infty}^0(x)$ in order to get a stationary process $(X_t)_{t\geq 0}$.  

\item
Assumption {\bf A2} is a conditional contraction property in $\L^p$ which is crucial for getting  the convergence of the backward iterations in $\L^p$ norms. Relaxing this assumption by introducing a random coefficient $\kappa_{t-1}$ instead of $\kappa$ can be problematic for getting a similar result. We discuss this point below.  

\end{enumerate}

\subsection{Comments on Assumption {\bf A2}}
Let us consider the dynamic (\ref{eq::discuss}), set $f_t(x)=F\left(x,Z_{t-1},\varepsilon_t\right)$, with $E=\R$, $p=1$ and assume that for every $t$, $\varepsilon_t$ is independent from $\mathcal{F}_{t-1}=\sigma\left((Z_j,\varepsilon_j): j\leq t-1\right)$. Then Assumption {\bf A2} is satisfied for $m=1$ if and only if there exists $\kappa\in (0,1)$ such that
$$\sup_{z\in E_1'}\E\left\vert F\left(x,z,\varepsilon_0\right)-F\left(y,z,\varepsilon_0\right)\right\vert\leq \kappa\vert x-y\vert.$$   
At a first sight, the latter condition is quite strong and it is natural to wonder if the following weaker assumption can be used, i.e. there exists a measurable function $\kappa:G\rightarrow (0,\infty)$ such that for every $z\in G$,
$$\E\left\vert F\left(x,z,\varepsilon_0\right)-F\left(y,z,\varepsilon_0\right)\right\vert\leq \kappa(z)\vert x-y\vert.$$
Of course, the challenging question concerns the convergence of the backward iterations when the function $\kappa$ may take values larger than $1$. 
However, a problem occurs for applying the successive contraction properties to the iterated random maps. Consider the iterations $f_t\circ f_{t-1}$. We have for $(x,y)\in E^2$,
\begin{eqnarray*}
\E\left\vert f_t\circ f_{t-1}(x)-f_t\circ f_{t-1}(y)\right\vert&=&\E\left[\E\left[\left\vert f_t\circ f_{t-1}(x)-f_t\circ f_{t-1}(y)\right \vert  \mathcal{F}_{t-1}\right]\right]\\
&\leq& \E\left[\kappa(Z_{t-1})\left\vert f_{t-1}(x)-f_{t-1}(y)\right\vert\right].
\end{eqnarray*}
If the random variable $\kappa(Z_{t-1})$ depends on past values of the error $\varepsilon_{t-j},j\geq 1$, 
it is stochastically dependent on the random map $f_{t-1}$ and also on $\mathcal{F}_{t-2}$. It is then not possible to use the contraction property of $f_{t-1}$ unless the function $\kappa$ can be bounded by a constant. To show that the successive iterations lose the memory with respect to initialization, this constant has to be smaller than $1$. Of course, this does not prove that the convergence of the iterations in $\L^1$ is not possible.

To show that the convergence of the backward iterations in $\L^1$ is problematic, we now consider a map $f_t$ linear in $x$,
a case for which an explicit solution is available. Accordingly, we assume that 
$$f_t(x)=\kappa(Z_{t-1})x+\varepsilon_t,$$
where the function $\kappa$ is bounded but not necessarily by $1$ and an integrable noise $\varepsilon_0$.
The dynamic is then given by an AR process with a random lag coefficient and it is widely known that the unique solution can be written as
\begin{equation}\label{dev}
X_t=\sum_{j\geq 1}\prod_{i=1}^j \kappa\left(Z_{t-i}\right)\varepsilon_{t-j}+\varepsilon_t,
\end{equation}
provided that $\E\log \kappa\left(Z_0\right)<0$. The series (\ref{dev}) converges almost surely. 
The $\L^1-$convergence is guaranteed from Theorem \ref{main}, as soon as $\kappa:=\Vert \kappa(Z_0)\Vert_{\infty}<1$, where for a random variable $X$, $\Vert X\Vert_{\infty}$ denotes its suppremum norm.
If $\kappa\geq 1$, convergence in $\L^1$ of the series \ref{dev} is much more difficult to get because of the possible stochastic dependence between the coordinates of the process $\left(Z_t\right)_{t\in\Z}$. 
Let us first note that such a problem occurs in the non ergodic case, when $Z_t=Z_0$ a.s. In this case, we have 
$$\E\left(\prod_{i=1}^j \kappa\left(Z_{t-i}\right)\vert \varepsilon_{t-j}\vert\right)=\E\left(\vert \varepsilon_0\vert\right)\cdot\E\left(\kappa(Z_0)^j\right)$$
and since $\E\left(\kappa(Z_0)^j\right)\geq \P\left(\kappa(Z_0)\geq 1\right)$, one cannot get convergence of the series (\ref{dev}) if 
$\P\left(\kappa(Z_0)\geq 1\right)>0$. In what follows, we also stress that a similar problem of convergence also occurs in the ergodic case. To this end, set $\phi(p)=\Vert \kappa(Z_0)\Vert_p$ for $p\geq 1$. The function $\phi$ is non decreasing and $\phi(\infty)=\Vert \kappa(Z_0)\Vert_{\infty}$.
Assumption {\bf A2} is satisfied as soon as $\Vert \kappa(Z_0)\Vert_{\infty}<1$.
It is then tempting to study the convergence of the solution only assuming that $\Vert \kappa(Z_0)\Vert_p< 1$ but $\Vert \kappa(Z_0)\Vert_q\geq 1$ for some $1\leq p<q$. However for any value of the pair $(p,q)$, there always exists an example of a process $\left(\kappa(Z_t)\right)_{t\in\Z}$ such that the series (\ref{dev}) is not converging in $\L^1$. To this end, assuming without loss of generality that $q$ is an integer, we define $\kappa(z)=z$ and $Z_{t-1}=a_{t-1}\cdots a_{t-q}$ where $(a_t)_{t\in\Z}$ is a process of i.i.d. nonnegative random variables, independent of $\left(\varepsilon_t\right)_{t\in\Z}$, and such that $\Vert a_0\Vert_p<1$ and $\Vert a_0\Vert_q\geq 1$. Since for $j\geq q$,
$$\prod_{i=1}^j \kappa\left(Z_{t-i}\right)=\prod_{i=1}^{q-1}a_{t-i}^i\cdot\prod_{i=q}^{j+1}a_{t-i}^q\prod_{i=0}^{q-2}a_{t-j-k+i}^{i+1},$$ 
we find 
$$\E\left(\prod_{i=1}^j \kappa\left(Z_{t-i}\right)\right)=\prod_{i=1}^{q-1}\E^2\left(a_0^i\right)\cdot \E^{j-q+2}\left(a_0^q\right).$$
Hence the previous expectation does not converge to $0$ when $j\rightarrow \infty$ and the series (\ref{dev}) cannot converge in $\L^1$.
The analysis of this linear case enlightens that a tail condition on $\kappa(Z_t)$ is not sufficient for getting this kind of convergence. In particular, the dependence structure of the process $\left(\kappa(Z_t)\right)_{t\in\Z}$ is also of major importance. This contrasts with AR processes with i.i.d. random coefficients, since in this case the condition $\E \kappa(Z_0)<1$ is necessary and sufficient for the convergence of the series (\ref{dev}) in $\L^1$. However, imposing an independence assumption on the covariate process is not reasonable.  

In the next section, we show that one can investigate a different mode of convergence for the iterations in model (\ref{eq::discuss}) and which allow to relax Assumption {\bf A2}. However, it is necessary to impose a strict exogeneity assumption on the covariate process.

\subsection{An additional result for strictly exogenous regressors}

In this subsection, we consider specifically equation (\ref{eq::discuss}) when the covariate process is independent of the error process.
In this case, conditionally on $Z$, the process is a time-inhomogeneous Markov chain. The terminology Markov chain in random environments 
is often used in the literature. See for instance \citet{Stenflo2}. 
The following result will not be central in the rest of the paper because substantial efforts could be needed to derive moment and weak dependence properties for the corresponding solution and it could be also difficult to obtain explicit conditions for dealing with  higher-order autoregressive processes.
This is why we only provide a result when the $f_t$'s satisfied a one-step contraction (i.e. $m=1$ in {\bf A2}). 
We assume that there exist a real number $p\geq 1$ and a state $x_0\in E$ such that the three following conditions are fulfilled. 

\paragraph{A0} 
The process $Z:=(Z_t)_{t\in\Z}$ is stationary and ergodic, $\left(\varepsilon_t\right)_{t\in\Z}$ is a process of i.i.d. random variables taking values in $E_2'$ and is independent of $Z$. 

\paragraph{A1'}
For every $z\in E'_1$, we have $\E\left[d^p\left(F(x_0,z,\varepsilon_1),x_0\right)\right]<\infty$.

\paragraph{A2'} 
There exists a measurable function $\kappa:E'\rightarrow (0,\infty)$ satisfying $\E\left(\log^{+}\kappa(Z_0)\right)<\infty$, $\E \log\kappa(Z_0)<0$ and such that for every $(x,y)\in E^2$,
$$\E\left[d^p\left(F(x,z,\varepsilon_0),F(y,z,\varepsilon_0)\right)\right]\leq \kappa^p(z)d^p(x,y).$$
Moreover,
$$\E\left[\log^{+}\int d^p\left(x_0,F(x_0,Z_0,u)\right)d\P_{\epsilon_1}(u)\right]<\infty.$$

We remind the notation $f_t(x)=F\left(x,Z_{t-1},\varepsilon_t\right)$. Here, we set $\mathcal{F}_t=\sigma\left((Z_j,\varepsilon_j):j\leq t\right)$ and $\E\left[X\vert Z\right]$ will denote the expectation of a random variable $X$ conditionally on the covariate process $(Z_t)_{t\in\Z}$. 
\begin{theo}\label{relax}
Suppose that Assumptions {\bf A0,A1'-A2'} hold. 
\begin{enumerate}
\item
For every $(x,t)\in E\times \Z$, there exists a random variable $X_t(x)$ such that $\E\left[d^p\left(X_t(x),x_0\right)\vert Z\right]<\infty$ a.s. and $\lim_{s\rightarrow \infty}\E\left[d^p\left(f_{t-s}^t(x),X_t(x)\right)\vert Z\right]=0$ a.s.
The sequence $\left(f^t_{t-s}(x)\right)_{s\geq 0}$ also converges almost surely to $X_t(x)$.
\item
For $x\neq y$, we have $\P\left(X_t(x)\neq X_t(y)\right)=0$. We then set $X_t=X_t(x)$. 
\item
The process $\left((X_t,Z_t)\right)_{t\in\Z}$ is stationary and ergodic.
\item
If $(Y_t)_{t\in\Z}$ is a non-anticipative process (i.e. $Y_t\in \mathcal{F}_t$) such that $\left((Y_t,Z_t)\right)_{t\in\Z}$ is stationary, for every $t\in\Z$, $Y_t=f_t\left(Y_{t-1}\right)$ and $\E\left[d^p\left(Y_0,x_0\right)\vert Z_0,Z_{-1},\ldots\right]<\infty$ a.s., then $Y_t=X_t$ a.s.
\end{enumerate}
\end{theo}

\paragraph{Notes}
\begin{enumerate}
\item
 The contraction inequality in Assumption {\bf A2'} can be restated as
$$\E\left[d^p\left(f_t(x),f_t(y)\right)\vert Z\right]\leq \kappa^p(Z_{t-1})d^p(x,y)\mbox{ a.s}.$$
It is then another example of contraction in conditional average. 
\item
In our context, our result can be seen as an improvement of Theorem $1$ given in \citet{Stenflo2} for Markov chains in random environments. In particular, we do not assume a uniform contraction with respect to the environment which is given by the exogenous process $Z$ in our random maps $f_t$. 
\end{enumerate}

\subsection{Example}\label{comp}
We compare the contraction conditions necessary to apply Theorem \ref{bibli}, Theorem \ref{main} or Theorem \ref{relax} on a specific example.
Let $\left((\varepsilon_t,Z_t)\right)_{t\in\Z}$ be a stationary sequence of pair of random variables and for $t\in\Z$, set $\mathcal{F}_t=\sigma\left((Z_j,\varepsilon_j):j\leq t\right)$. Assume that 
$\E\left(\varepsilon_1^2\vert \mathcal{F}_0\right)=1$ and $\E\left(\varepsilon_1\vert\mathcal{F}_0\right)=0$. For $i=1,2$, let $a_i:E_1'\rightarrow\R$ and $b_i:E_1'\rightarrow \R_+$ be some measurable maps such that $\E\log^{+}a_i(Z_0)<\infty$ and $\E\log^{+}b_i(Z_0)<\infty$. We consider the following AR-ARCH model with functional coefficients 
$$X_t=a_0(Z_{t-1})+a_1(Z_{t-1})X_{t-1}+\varepsilon_t\sqrt{b_0(Z_{t-1})+b_1(Z_{t-1})X_{t-1}^2}.$$
Here we set $E=\R$ and $d(x,y)=\vert x-y\vert$ and $p=2$. 
Setting $\sigma(z,x)=\sqrt{b_0(z)+b_1(z)x^2}$, one can note that $\left\vert \sigma(z,x)-\sigma(z,y)\right\vert\leq \sqrt{b_1(z)}\vert x-y\vert$.
\begin{enumerate}
\item
To apply Theorem \ref{bibli}, we compute the Lipschitz constant $c(f_1)$ of the random map $f_1$. We have 
$$c(f_1)=\sup_{v\in\R}\vert f_1'(v)\vert=\sup_{v\in \R}\left\vert a_1(Z_0)+\frac{\varepsilon_1b_1(Z_0)v}{\sqrt{b_0(Z_0)+b_1(Z_0)v^2}}\right\vert.$$
Making the change of variable $\overline{v}=sign(\varepsilon_1) sign\left(a_1(Z_0)\right)v$, we have
$$c(f_1)=\sup_{\overline{v}\in\R}\left[\left\vert a_1(Z_0)\right\vert+\frac{\left\vert \varepsilon_1\right\vert b_1(Z_0)\overline{v}}{\sqrt{b_0(Z_0)+b_1(Z_0)\overline{v}^2}}\right\vert.$$
We then obtain $c(f_1)=\left\vert a_1(Z_0)\right\vert+ \sqrt{b_1(Z_0)}\vert\varepsilon_1\vert$ and
Theorem \ref{bibli} applies as soon as 
\begin{equation}\label{bench1}
\E \log c(f_1)<0.
\end{equation}
\item
To apply Theorem \ref{main}, note that
\begin{eqnarray*}
\E\left[\left\vert f_t(x)-f_t(y)\right\vert^2 \vert \mathcal{F}_{t-1}\right]&=&a_1(Z_{t-1})^2(x-y)^2+\left(\sigma\left(Z_{t-1},x\right)-\sigma\left(Z_{t-1},x\right)\right)^2\\
&\leq& \left(a_1(Z_{t-1})^2+b_1(Z_{t-1})\right)\cdot\vert x-y\vert^2.
\end{eqnarray*}
One can then show that Theorem \ref{main} applies as soon as 
\begin{equation}\label{bench2}
\kappa^2:=\sup_z\left(a_1(z)^2+b_1(z)\right)<1. 
\end{equation}
\item
If the two processes $(\varepsilon_t)_{t\in\Z}$ and $(Z_t)_{t\in\Z}$ are independent and the $\varepsilon_t'$s are i.i.d., Theorem \ref{relax} applies as soon as 
\begin{equation}\label{bench3}
\E\log\left(a_1(Z_0)^2+b_1(Z_0)\right)<0. 
\end{equation}
\end{enumerate}
Under the strict exogeneity assumption, we note that (\ref{bench3}) is weaker than (\ref{bench2}). 
However, (\ref{bench2}) ensures the existence of a second order moment for the solution whereas (\ref{bench3}) only guarantees that $\E\left(X_t^2\vert Z\right)<\infty$ a.s.
On the other hand, (\ref{bench1}), which only ensures existence of a stationary solution, is not necessarily weaker than (\ref{bench2}) or (\ref{bench3}). For instance, if the noise process has a Rademacher distribution, $\P(\varepsilon_1=1)=\P(\varepsilon_1=-1)=1/2$ and 
the functional coefficients are constant, (\ref{bench1}) writes $\vert a_1\vert +\sqrt{b_1}<1$ which is more restrictive than (\ref{bench2}) or (\ref{bench3}). But if $a_1$ is identically equal to $0$, (\ref{bench1}) writes as $\frac{1}{2}\E\log b_1(Z_0)+\E\log \vert \varepsilon_0\vert<0$ which is weaker than (\ref{bench3}), since from Jensen's inequality, we have
$$\E\log\vert\varepsilon_0\vert\leq \log \E\vert \varepsilon_0\vert\leq \log\E^{1/2}\left(\varepsilon_0^2\right)=0.$$

\subsection{A result for higher-order autoregressions}\label{HO}
In this subsection, we extend Theorem \ref{main} to higher-order autoregressive processes. We only consider stationary processes in this part. The main result, Theorem \ref{autoreg}, is particularly interesting for multivariate autoregressions for which Lipschitz type properties can be obtained equation by equation. See Section \ref{examples} for an application of Theorem \ref{autoreg} to various examples.  

For a given real number $0<o\leq 1$, we define the distance $\Delta$ on $\R$, by $\Delta(u,v)=\vert u-v\vert^o$ for $u,v\in\R$. 
Let $E$ be a subset of $\R^k$ and $\Vert\cdot\Vert$ an arbitrary norm on $\R^k$.  
Our aim is to study existence of solutions for 
the following recursive equations:
\begin{equation}\label{recur}
X_t=F\left(X_{t-1},\ldots,X_{t-q},\zeta_t\right),\quad t\in\Z,
\end{equation}
where $F:E^q \times E'\rightarrow E$ is a measurable function. 
Note that one can always associate a random map $f_t$ on $E^q$ to the dynamic (\ref{recur}). To this end, for $t\in \Z$ and $x_1,\ldots,x_q\in E$, we set 
$$f_t(x_1,\ldots,x_q)=\left(F\left(x_1,\ldots,x_q,\zeta_t\right),x_1,\ldots,x_{q-1}\right).$$

We first introduce additional notations. We denote by $\mathcal{M}_k$ the set of square matrices with real coefficients and $k$ rows and if $A\in\mathcal{M}_k$, $\rho(A)$ the spectral radius of the matrix $A$. Moreover,
for $x,y\in \R^k$ and $p\geq 1$, the vector $\left(\Delta^p(x_1,y_1),\ldots,\Delta^p(x_k,y_k)\right)'$ will be denoted by $\Delta^p_{vec}(x,y)$. Finally, we introduce a partial order relation $\preceq$ on $\R^k$ and such that $x\preceq y$ means $x_i\leq y_i$ for $i=1,\ldots,k$.

The following assumptions will be needed.
\begin{description}
\item[B1]
The process $\left(\zeta_t\right)_{t\in\Z}$ is stationary and ergodic adapted to a filtration $\left(\mathcal{F}_t\right)_{t\in\Z}$.
\item [B2]
For any $y\in E^q$, $\E\left[\Vert F\left(y,\zeta_1\right)\Vert^{op}\right]<\infty$.
\item [B3] 
There exist some matrices $A_1,\ldots,A_q\in \mathcal{M}_k$ with nonnegative elements, satisfying $\rho\left(A_1+\cdot+A_q\right)<1$ and  such that for $y,y'\in E^q$,
$$\E\left[\Delta^p_{vec}\left(F(y,\zeta_t),F(y',\zeta_t)\right)\vert\mathcal{F}_{t-1}\right]\preceq \sum_{i=1}^q A_i \Delta^p_{vec}(y_i,y'_i).$$
\end{description}

Though the following result is stated for an arbitrary pair $o\in (0,1),p\geq 1$, the two interesting cases are $o\in (0,1), p=1$ and $o=1,p>1$.
\begin{theo}\label{autoreg}
Suppose that Assumptions {\bf B1-B3} hold. There then exists a unique stationary and non-anticipative process $(X_t)_{t\in\Z}$ solution of (\ref{recur}) and such that $\E\left[\Vert X_t\Vert^{op}\right]<\infty$. Moreover, this process is ergodic.
\end{theo}

\paragraph{Notes} 
\begin{enumerate}
\item
If $\zeta_t=(Z_{t-1},\varepsilon_t)$ with $\varepsilon_t$ independent of $\mathcal{F}_{t-1}=\sigma\left((Z_j,\varepsilon_j):j\leq t\right)$ and $k=o=1$, Assumption {\bf B3} writes 
$$\sup_z\E\left[\left\vert F(y_1,\ldots,y_q,z,\varepsilon_1)-F(y_1',\ldots,y_q',z,\varepsilon_1)\right\vert^p\right]\leq \sum_{i=1}^q A_i\vert y_i-y'_i\vert^p,$$
with $\rho\left(\sum_{i=1}^q A_i\right)=\sum_{i=1}^q A_i<1$. This provides a quite simple criterion for application to autoregressive processes.
\item
In the spirit of Section \ref{comp}, the previous criterion can be checked for models with varying parameters, directly constructed from smooth parametric autoregressive processes. Consider the model 
$$Y_t=\overline{F}_{\theta}\left(Y_{t-1},\ldots,Y_{t-q},\varepsilon_t\right),\quad t\in\Z, \theta\in\Theta\subset \R^{\overline{e}}.$$
If 
$$\E\left[\left\vert \overline{F}_{\theta}(y_1,\ldots,y_q,\varepsilon_1)-\overline{F}_{\theta}(y_1',\ldots,y_q',\varepsilon_1)\right\vert^p\right]\leq \sum_{i=1}^q \overline{A}_i(\theta)\vert y_i-y'_i\vert^p,$$
the model
$$Y_t=\overline{F}_{\theta(Z_{t-1})}\left(Y_{t-1},\ldots,Y_{t-q},\varepsilon_t\right)$$
satisfies {\bf B3} as soon as $\sum_{i=1}^q \Vert \overline{A}_i\Vert_{\infty}<1$. 
We then obtain a model with exogenous covariates by replacing parameter $\theta$ with a varying parameter $\theta\left(Z_{t-1}\right)$ where $\theta:\R^e\rightarrow \Theta$ is a measurable map. See also the note after Proposition \ref{CHARME} for a discussion.  
\end{enumerate}

\section{Functional dependence measure}\label{mix}

The functional dependence measure has been introduced by \citet{Wu} and is particularly interesting for autoregressive processes
which are not necessarily strong mixing or for which getting strong mixing conditions requires additional regularity conditions on the noise distribution. The single requirement is to get a Bernoulli shift representation of the stochastic process of interest, i.e. $X_t=H\left(\xi_t,\xi_{t-1},\ldots\right)$ where $(\xi_t)_{t\in\Z}$ is a sequence of i.i.d. random variables taking values in a measurable space $\left(G,\mathcal{G}\right)$. The functional dependence measure is expressed in terms of some coefficients 
which evaluate for $t\geq 0$ the $\L^p-$distance between $X_t$ and a copy $\overline{X}_t$, obtained by replacing $\xi_0$ with $\xi_0'$,  $\xi'_0$ following the same distribution as $\xi_0$ and being independent from the sequence $\left(\xi_t\right)_{t\in\Z}$. See below for the definition of these coefficients. 
Most of the limit theorems and deviation inequalities have been derived under such dependence measures. See for instance \citet{WuStrong} and \citet{WuWu}. 
Such asymptotic results have been applied to various statistical problems. See for instance \citet{WuH} for kernel estimation for time series, \citet{WuHan} for covariance estimation or \citet{LW} for spectral density estimation.    
The notion of functional dependence is then an attractive alternative to the usual strong mixing when the process is defined by stochastic recursions.
Our aim in this section is to show that under the assumptions of Theorem \ref{main} or Theorem \ref{autoreg}, when the process $\left(\zeta_t\right)_{t\in\Z}$ or more specifically the covariate process $(Z_t)_{t\in\Z}$ satisfies this kind of dependence, the functional dependence measure of the solution $(X_t)_{t\in\Z}$ can be controlled. We will then provide a new wide class of examples for which the aforementioned references provide an important number of statistical applications. 
In this section, we assume that the state space $E$ is a subspace of $\R^k$ and the distance $d$ is given by 
$d(x,y)=\vert x-y\vert^{o}$ where $\vert\cdot\vert$ is an arbitrary norm on $\R^k$ and $0<o\leq 1$.

\subsection{Dependence coefficients for general iterated random functions}

Assume that the process $\zeta$ has a Bernoulli shift representation, i.e. $\zeta_t=H\left(\xi_t,\xi_{t-1},\ldots\right)$ for some measurable map $H:G^{\N}\rightarrow E'$ and $\left(\xi_t\right)_{t\in\Z}$ is a sequence of i.i.d., $G-$valued random variables. We then have for every $(x,t)\in E\times \Z$,
$$f_t(x)=F\left(x,H\left(\xi_t,\xi_{t-1},\ldots\right)\right)$$
and the map $f_t$ has itself a Bernoulli shift representation. 
To define the functional measure coefficients, we then define a new sequence $\left(\overline{\xi}_t\right)_{t\in\Z}$ such that $\overline{\xi}_0=\xi_0'$ and $\overline{\xi}_t=\xi_t$ for $t\neq 0$. Here $\xi_0'$ is a copy of $\xi_0$ which is assumed to be independent from $(\xi_t)_{t\in\Z}$.
 Moreover, for $t>0$, let
$$\overline{f}_t(x)=F(x,H\left(\overline{\xi}_t,\overline{\xi}_{t-1},\ldots)\right).$$
we define for $t\geq 0$ and $p\geq 1$, 
$$\theta_{p,t}=\E^{1/p}\left[d\left(f_{-\infty}^t,\overline{f}_{-\infty}^t\right)^p\right].$$
Moreover, for $h\in\N$, let $\Theta_{p,h}=\sum_{t\geq h}\theta_{p,t}$. Two cases of interest are $p=1$ and $o=1,p>1$.

Our aim is to get an upper bound for the functional dependence coefficients $\Theta_{p,h}$. To this end, we add other assumptions.
Here we set for $t\in \Z$, $\mathcal{F}_t=\sigma\left(\xi_{t-j}:j\geq 0\right)$.

\paragraph{A3} There exists a measurable function $S:E\rightarrow \R_+$ and $r,s\geq p$ such that $r^{-1}+s^{-1}=p^{-1}$ and $\E\left[S(f_{-\infty}^0)^s\right]<\infty$ and 
for all $x\in E$ and $t\geq 1$,
$$\E^{1/p}\left[d^p\left(\overline{f}_t(x),f_t(x)\right)^p\vert \sigma(\xi_0')\vee\mathcal{F}_{t-1}\right]\leq S(x)H_{t-1},$$
where $H_{t-1}$ is a random variable measurable with respect to $\sigma(\xi_0')\vee\mathcal{F}_{t-1}$ and such that $\E\vert H_{t-1}\vert^r<\infty$. 
 
An immediate consequence of Assumption {\bf A4} is that for any random variable $V_{t-1}$, measurable with respect to $\sigma(\xi_0')\vee \mathcal{F}_{t-1}$, we have $\E\left[d^p\left(\overline{f}_t(V_{t-1}),f_t(V_{t-1})\right)\right]\leq \E\left[S(V_{t-1})^pH_{t-1}^p\right]$ and from H\"older's inequality, we get 
$$\Vert d\left(\overline{f}_t(V_{t-1}),f_t(V_{t-1})\right)\Vert_p\leq \Vert S(V_{t-1})\Vert_s \Vert H_{t-1}\Vert_r.$$

When $\zeta_t=\left(Z_{t-1},\varepsilon_t\right)$, with $Z_t$ taking values in a Borel subset $E_1'$ of $\R^e$ and $\varepsilon_t$ taking values in a Polish space $E_2'$, we will still denote by $\vert\cdot\vert$ an arbitrary norm on $\R^e$ and we also set $d(z,z')=\vert z-z'\vert^o$ for $z,z'\in\R^e$ to avoid additional notations. We will use two specific assumptions.

\paragraph{A3'}
There exists a measurable function $S:E\rightarrow \R_+$ and $r,s\geq p$ such that $r^{-1}+s^{-1}=p^{-1}$ and $\E\left[S(f_{-\infty}^0)^s\right]<\infty$, $\E\vert Z_0\vert^{ro}<\infty$ and 
for all $x\in E$ and $t\geq 1$,
$$\E^{1/p}\left[d^p\left(F\left(x,z,\varepsilon_0\right),F\left(x,z',\varepsilon_0\right)\right)\right]\leq S(x)d(z,z'),$$

\paragraph{A4}
Let $(\eta_t)_{t\in\Z}$ be a sequence of random variables taking values in a measurable space $(G_1,\mathcal{G}_1)$ and such that $Z_t=H'\left(\eta_t,\eta_{t-1},\ldots\right)$ for a measurable function $H'$. Moreover, setting $\xi_t=(\varepsilon_t,\eta_t)$, we assume that $\left(\xi_t\right)_{t\in\Z}$ is a sequence of i.i.d. of random variables taking values in $G=E_2'\times G_1$.

From {\bf A4}, we have the specific representation
$$\zeta_t=H\left(\xi_t,\xi_{t-1},\ldots\right):=\left(H'(\eta_{t-1},\eta_{t-2},\ldots),\varepsilon_t\right),\quad t\in\Z.$$
The map $H$ takes values in $E'=E_1'\times E_2'$. Note that our formulation allows the covariate process to have a general form,
including a VARMA or GARCH process among others.

\begin{prop}\label{WD}
\begin{enumerate}
\item
Suppose that Assumptions {\bf A1-A4} hold. 
For any $h\geq 2$, there then exists $C_1>0$ not depending on $h$ such that 
\begin{equation}\label{bbound}
\Theta_{p,h}\leq  C_1\left[\kappa^{h/m}+\sum_{i=0}^{h-1}\kappa^{i/m}\eta_{r,h-i}+\sum_{i\geq h}\kappa^{i/m}\eta_{r,1}\right],
\end{equation}
with 
$$\eta_{r,j}=\sum_{t\geq j}\Vert H_{t-1}\Vert_r,\quad j\geq 1.$$
In particular, if $\eta_{r,1}<\infty$, there exists $C>0$, not depending on $h$, such that 
\begin{equation}\label{goodmaj}
\Theta_{p,h}\leq C\left[\kappa^{h/m}+\sum_{i=0}^{h-1}\kappa^{i/m}\eta_{r,h-i}\right].
\end{equation}
\item
Suppose that Assumptions {\bf A1-A2-A4} and {\bf A3'} hold.
Then we get the bound (\ref{goodmaj}) with $\eta_{r,j}=\Theta_{r,j-1}(Z)$.
\end{enumerate}
\end{prop}

\paragraph{Notes} 
\begin{enumerate}
\item
Let us comment on Assumption {\bf A4}. Under this assumption, the $\eta_t'$s are i.i.d. as well as the $\varepsilon_t'$s and for any $t\in\Z$, $\varepsilon_t$ is independent from $\mathcal{F}_{t-1}=\sigma\left(\xi_s:s\leq t-1\right\}$. Note that we allow simultaneous dependence between $\varepsilon_t$ and $\eta_t$. For instance, we can set $\eta_t=K(\varepsilon_t,U_t)$ where $K$ is a measurable function and $U$ is a sequence of i.i.d. random variables, independent from the sequence $\varepsilon$.
This assumption is then more flexible than the complete independence between the two error processes $\varepsilon$ and $\eta$, which implies strict exogeneity. 
\item
It can happen that our assumptions are satisfied with some $p$, leading to an upper bound for the functional dependence coefficients $\theta_{p,t}$, while it is required a condition on $\theta_{q,t}$ or $\Theta_{q,h}$ for $q>p$ for applying some limit theorems or statistical results.
This is still possible if one can prove finiteness of higher-order moments for the solution, e.g. if $\E\left[\vert f_{-\infty}^0\vert^{q'o}\right]<\infty$ for some $q'>q$. Indeed, from H\"{o}lder's inequality, we have
$$\theta_{q,t}\leq \theta_{p,t}^{\frac{p(q'-q)}{q(q'-p)}}\theta_{q',t}^{\frac{q'(q-p)}{q(q'-p)}}.$$
Moreover, $\theta_{q',t}\leq 2\Vert \left\vert f_{-\infty}^0\right\vert^o\Vert_{q'}$.

\end{enumerate}

\subsection{Dependence coefficients for higher-order autoregressions with exogenous covariates}

Here, we revert to higher-order autoregressions considered in Section \ref{HO}. 
We consider directly the case $\zeta_t=\left(Z_{t-1},\varepsilon_t\right)$ with Assumption {\bf A4} being satisfied. Additionally to Assumptions {\bf B1-B2-B3} and {\bf A4}, the following assumption, which is the analogue of {\bf A3'}, will be needed. 

\begin{description}
\item[B4] 
If $r>0$ and $s\in \R_+\cup\{\infty\}$ are such that $r^{-1}+s^{-1}=p^{-1}$, there exists a measurable function $S:E^q\rightarrow \R$ such that $\E\left[S(X_q,\ldots,X_1)^s\right]<\infty$, $\E\vert Z_0\vert^{or}<\infty$ and for all $z,z'\in E_1'$ and $x_1,\ldots,x_q\in E^q$,
$$\E^{1/p}\left[d^p\left(F(x_1,\ldots,x_q,z,\varepsilon_0),F(x_1,\ldots,x_q,z',\varepsilon_0)\right)\right]\leq S(x_1,\ldots,x_q)d(z,z').$$
\end{description}

The result is the following.

\begin{prop}\label{FDM}
Suppose that Assumptions ${\bf B2-B4}$ and ${\bf A4}$ hold. There then exists $C>0$ and $\rho\in(0,1)$, such that for all $h\geq 1$, $\Theta_{p,h}(X)\leq C \left[\rho^h+\sum_{i=1}^h\rho^i\Theta_{r,h-i}(Z)\right]$.
\end{prop}

\paragraph{Note.}
From the upper bound given in Proposition \ref{FDM}, we note that the decay of $\Theta_{p,h}(X)$ is polynomial (respectively geometric) in $h$ if the decay of $\Theta_{r,h}$ is polynomial (respectively geometric) in $h$.

\subsection{A central limit theorem}\label{TLC}
To illustrate the usefulness of our results, we give below a central limit theorem for partial sums
$$S_n=\sum_{i=1}^nf\left(X_t,X_{t-1},\ldots,X_{t-k}\right),$$
where $f$ is some real-valued measurable function and $(X_t)_{t\in\Z}$ is a stochastic process solution of $X_t=F\left(X_{t-1},\ldots,X_{t-q},\zeta_t\right)$ and the assumptions of either Proposition \ref{WD} or Proposition \ref{FDM} are satisfied.
The following result, which is a straightforward corollary of the invariance principle given in \citet{Wu}, is not the most general as possible. In particular, when a moment of order greater than $p$ is available for the stationary solution, different assumptions on the function $f$ could be used.

\begin{theo}\label{TLC}
Suppose that either Assumptions {\bf A1-A4} or Assumptions {\bf B2-B3-B4-A4} hold true for some $p>2$ and $\Theta_{p,0}(X)<\infty$. If there exists $C>0$ and $0\leq \ell\leq \frac{p-2}{2}$ such that for $x_i,x_i'\in E$, $0\leq i\leq k$,
$$\left\vert f\left(x_0,\ldots,x_k\right)-f\left(x_0',\ldots,x_k'\right)\right\vert\leq C\left[1+\sum_{i=0}^k\left(\vert x_i\vert^{o\ell}+\vert x_i'\vert^{o\ell}\right)\right]\cdot\sum_{i=0}^kd(x_i,x_i').$$
Then we have the weak convergence
$$\frac{1}{\sqrt{n}}\left(S_n-\E S_n\right)\Rightarrow \mathcal{N}\left(0,\sigma^2\right),$$
with $\sigma^2=\sum_{j\in\Z}\c\left(Y_j,Y_0\right)$, $Y_t=f\left(X_t,\ldots,X_{t-k}\right)$.
\end{theo}

\section{Examples}\label{examples}

\subsection{CHARN models}
In this section, we consider conditional heteroscedastic autoregressive nonlinear (CHARN) models
such as in \citet{Tsyb1} or \citet{Tsyb2} but that can encompass exogenous regressors. More precisely, 
we consider the dynamic
\begin{equation}\label{CH}
Y_t=f_1\left(Y_{t-1},\ldots,Y_{t-q},Z_{t-1}\right)+\varepsilon_t f_2\left(Y_{t-1},\ldots,Y_{t-q},Z_{t-1}\right),
\end{equation}
where $q$ is a positive integer and $f_1,f_2:\R^q\times E_1'\rightarrow \R$ are measurable functions.
In order to study stationary solutions of the recursions (\ref{CH}), the following assumptions will be needed.

\begin{description}
\item[CH1]
The process $\left((Z_t,\varepsilon_t)\right)_{t\in\Z}$ is stationary and ergodic.
\item[CH2]
For $j=1,2$, there exist measurable functions $a_{i,j}:\R^d\times \R_+$, $1\leq i\leq q$ such that 
$$\left\vert f_j(y_1,\ldots,y_q,z)-f_j(y'_1,\ldots,y'_q,z)\right\vert\leq \sum_{i=1}^q a_{i,j}(z)\left\vert y_i-y'_i\right\vert.$$
\item[CH3]
There exist a real number $p\geq 1$ such that $\Vert\varepsilon_1\Vert_p<\infty$ and $r,s\geq p$ such that $r^{-1}+s^{-1}=p^{-1}$, $s$ can be infinite, $\E\vert Z_0\vert^r<\infty$ and
two functions $L_1,L_2$ defined on $\R^q$ and such that for $j=1,2$, $y_1,\ldots,y_q\in\R$ and $z,z'\in E_1'$, 
$$\left\vert f_j(y_1,\ldots,y_q,z)-f_j(y_1,\ldots,y_q,z')\right\vert\leq L_j(y_1,\ldots,y_q)\vert z-z'\vert.$$
\end{description}

For $t\in\Z$ and $i=1,\ldots,q$, we set $c_{i,t}=a_{i,1}\left(Z_{t-1}\right)+a_{i,2}\left(Z_{t-1}\right)\vert\varepsilon_t\vert$,
We then define a sequence of random matrices ${\bf A}=(A_t)_{t\in\Z}$ by 
$$A_t=\begin{pmatrix} c_{1,t}& c_{2,t}&\cdots & c_{q,t}\\
&&&0\\
&I_{q-1}&&\vdots\\
&&&0
\end{pmatrix}.$$
Finally, we denote by $\chi({\bf A})$ the Lyapunov exponent of the sequence ${\bf A}$, i.e. 
$$\chi({\bf A})=\lim_{n\rightarrow \infty}\frac{\E\left[\log\Vert A_n\cdots A_1\Vert\right]}{n},$$
where $\Vert\cdot\Vert$ is an arbitrary norm on the space of square matrices of size $q\times q$.

\begin{prop}\label{CHARME}
Suppose that Assumptions {\bf CH1-CH2} hold. 
\begin{enumerate}
\item
Suppose that $\chi({\bf A})<0$. There then exists a unique stationary process $(Y_t)_{t\in\Z}$ solution of (\ref{CH}) which is also ergodic.

\item
Assume additionally that for every $t\in\Z$, $\varepsilon_t$ is independent from $\mathcal{F}_{t-1}$. If there exist $\overline{x}\in\R^q$, $p\geq 1$ such that $f_1(\overline{x},Z_0)+\varepsilon_1 f_2(\overline{x},Z_0)\in \L^p$ and 
\begin{equation}\label{newdeal}
\sum_{i=1}^q \sup_z\Vert a_{i,1}(z)+a_{i,2}(z)\vert \varepsilon_1\vert \Vert_p<1,
\end{equation}
there then exists a unique stationary and non-anticipative process solution of (\ref{CH}) which is also ergodic and such that $\E\vert Y_1\vert^p<\infty$.

\item
Assume furthermore that Assumptions ({\bf CH3}) and ({\bf A4}) hold true with $\E L^s_j(Y_q,\ldots,Y_1)<\infty$ for $j=1,2$. There then exists $C>0$ and $\rho\in(0,1)$, such that for all $h\geq 1$, $\Theta_{p,h}(X)\leq C \left[\rho^h+\sum_{i=1}^h\rho^i\Theta_{r,h-i}(Z)\right]$.
\end{enumerate}
\end{prop}

\paragraph{Notes} 
\begin{enumerate}
\item
Our results can be useful for dealing with models with functional coefficients in the spirit of the example given in Section \ref{comp}.
See also the notes after the statement of Theorem \ref{autoreg}.
For $j=1,2$, let $m^{(j)}_{\theta}:\R^q\rightarrow\R$ be some functions depending on some parameters $\theta \in \R^e$ and such that
$$\left\vert m^{(j)}_{\theta}(y_1,\ldots,y_q)-m^{(j)}_{\theta}(y'_1,\ldots,y'_q)\right\vert\leq \sum_{i=1}^q d_{i,j}(\theta)\left\vert y_i-y'_i\right\vert,$$
for some nonnegative real numbers $d_{i,j}(\theta)$, $1\leq i\leq q$. If $\theta$ is replaced by a function $\theta(\cdot):\R^d\rightarrow \R^e$ and 
$f_j(y_1,\ldots,y_q,z)=m^{(j)}_{\theta(z)}(y_1,\ldots,y_q)$ for $j=1,2$,
one can then consider some standard autoregressive processes and obtain a version with functional parameters depending on exogenous covariates. For instance, threshold autoregressions or power-ARCH volatility, 
$$m^{(1)}_{\theta}(y_1,\ldots,y_q)=\theta_0+\sum_{i=1}^q\left(\theta_i y_i^{+}+\theta_{i+q}y_i^{-}\right),\quad 
m^{(2)}_{\theta}(y_1,\ldots,y_q)=\left(\theta_0+\sum_{i=1}^q \theta_i\vert y_i\vert^{\delta}\right)^{1/\delta},$$
where $x^+$ and $x^{-}$ denotes respectively the positive part and the negative part of a real number $x$
and $\delta\geq 1$.

\item
When $\varepsilon_t$ is not necessarily independent from $\mathcal{F}_{t-1}$, (\ref{newdeal}) can be replaced with the following more abstract condition. There exists $\eta\in (0,1)$ such that 
\begin{equation}\label{newdeal2}
\sum_{i=1}^q \E^{1/p}\left[ \left(a_{i,1}(Z_0)+a_{i,2}(Z_0)\vert\varepsilon_1\vert\right)^p\vert\mathcal{F}_0\right]\leq 1-\eta\mbox{ a.s.}
\end{equation}
For $q=1$, let us compare (\ref{newdeal2}) with the condition $\chi({\bf A})<0$, which reduces to $$\E\left[\log\left(a_{1,1}(Z_0)+a_{1,2}(Z_0)\vert\varepsilon_1\vert\right) \right]<0.$$ 
This latter condition is much weaker than (\ref{newdeal2}). 
Indeed, (\ref{newdeal2}) entails that 
$$\E\left[a_{1,1}(Z_0)+a_{1,2}(Z_0)\vert\varepsilon_1\vert\right]\leq \Vert a_{1,1}(Z_0)+a_{1,2}(Z_0)\vert\varepsilon_1\vert\Vert_p<1$$ 
and from Jensen's inequality, $\log\E\left[a_{1,1}(Z_0)+a_{1,2}(Z_0)\vert\varepsilon_1\vert\right]\leq \chi({\bf A})$.

For $q\geq 2$, it is more difficult to obtain explicit conditions which guaranty that $\chi({\bf A})<0$.

\item
Using the results of \citet{WuLiu}, a nonparametric kernel estimation of the functions $f$ and $g$ is possible. Proposition \ref{CHARME} gives precise assumptions under which it is possible to control the functional dependence measure of some CHARN models when the regressors include lag values of the response as well as exogenous covariates. We then obtain additional examples of time series models for which standard nonparametric estimators of the regression function are still consistent.
\end{enumerate}

\subsection{GARCH processes}
GARCH processes with exogenous regressors have been considered recently by \citet{GG} or \citet{Francq}.
We consider here the asymmetric power GARCH studied by \citet{Francq}. 
The model is defined as follows.
\begin{equation}\label{defGARCH}
Y_t=\varepsilon_t h_t^{1/\delta},\quad
h_t=\pi' Z_{t-1}+\sum_{i=1}^q\left\{\beta_i h_{t-i}+\alpha_{i+}(Y_{t-i}^{+})^{\delta}+\alpha_{i-}(Y_{t-i}^{-})^{\delta}\right\},
\end{equation}
where $(\varepsilon_t)_{t\in\Z}$ and $(Z_t)_{t\in\Z}$ are two sequences of random variables taking 
values in $\R$ and $\R_+^d$ respectively, $\delta>0$, $\pi\in\R_+^d$ and the $\beta_i'$s, $\alpha_{i+}'$s and $\alpha_{i-}'$s are nonnegative real numbers.
Optimal stationarity properties of time series models defined by (\ref{defGARCH})
have been obtained by \citet{Francq}, using a version of Theorem \ref{bibli} for affine random maps. 
In contrast, we use our results to get existence of a moment of order $\delta$ for the unique stationary solution.
The following assumptions will be needed.

\begin{description}
\item[G1] The process $\left((Z_t,\varepsilon_t)\right)_{t\in\Z}$ is stationary and ergodic and $\E\vert Z_0\vert<\infty$.
\item [G2] There exist $s_{-},s_{+}$ such that $\E\left[(\varepsilon_t^{+})^{\delta}\vert \mathcal{F}_{t-1}\right]\leq s_{+}$ and $\E\left[(\varepsilon_t^{-})^{\delta}\vert \mathcal{F}_{t-1}\right]\leq s_{-}$ a.s.
and $\gamma:=\sum_{i=1}^p\left(\beta_i+s_{+}\alpha_{i+}+s_{-}\alpha_{i-}\right)<1$.
\end{description}

\begin{prop}\label{Garch}
Suppose that Assumptions {\bf G1-G2} hold. 
\begin{enumerate}
\item
There then exists a unique stationary and non-anticipative 
solution $(Y_t)_{t\in\Z}$ of (\ref{defGARCH}). This solution is ergodic and satisfies $\E\vert Y_0\vert^{\delta}<\infty$. 
\item
Additionally, assume that Assumption {\bf A4} holds true. Let $H_t=\left((Y_t^{+})^{\delta},(Y_t^{+})^{\delta},h_t\right)$. There then exists $C>0$ and $\rho\in(0,1)$, such that for all $h\geq 1$, $\Theta_{1,h}(H)\leq C \left[\rho^h+\sum_{i=1}^h\rho^i\Theta_{1,h-i}(Z)\right]$. Moreover, if $\delta\geq 1$, we have the bound 
$$\theta_{\delta,t}(Y)\leq \theta_{1,t}^{1/\delta}\left(Y^{+}\right)+\theta_{1,t}^{1/\delta}\left(Y^{-}\right),\quad t\in \N.$$
\end{enumerate}
\end{prop}

\paragraph{Note.} Let us consider the example of a GARCH process. We then set $\delta=2$, $\alpha_{j+}=\alpha_{j-}=\alpha_j$ and we assume that $\left(\varepsilon_t\right)_{t\in\Z}$ is a martingale difference, adapted to the filtration $\left(\mathcal{F}_t\right)_{t\in\Z}$ with $\mathcal{F}_t=\sigma\left((\varepsilon_s,Z_s): s\leq t\right)$. 
Set $v_{t-1}=\E\left[\varepsilon_t^2\vert\mathcal{F}_{t-1}\right]$. If there exists a positive real number $v_+$ such that $v_{t-1}\leq v_+$ a.s., the contraction condition in {\bf G2} reduces to
$v_+\sum_{j=1}^q\left(\alpha_j+\beta_j\right)<1$. For GARCH processes with i.i.d. innovations $\varepsilon_t$, we recover a standard condition ensuring the existence of a solution with a finite second moment.

\subsection{Poisson autoregressions}\label{Poissss}
We consider the PARX model introduced in \citet{Cav}. The idea is to model 
the conditional distribution of $Y_t$ given $\mathcal{F}_{t-1}$ by a Poisson distribution with a random intensity $\lambda_t$ depending on past values and a covariate process.
More precisely, we assume that
\begin{equation}\label{Poisss}
Y_t=N^{(t)}_{\lambda_t},\quad \lambda_t=\beta_0+\sum_{j=1}^q\beta_j \lambda_{t-j}+\sum_{j=1}^q\alpha_jY_{t-j}+\pi'Z_{t-1},
\end{equation}
where $\left(N^{(t)}\right)_{t\in\Z}$ is a sequence of i.i.d. Poisson processes with intensity $1$, $\beta_0,\ldots,\beta_q$, $\alpha_1,\ldots,\alpha_q$ are nonnegative real numbers and $\pi$ is a vector of $\R^d$ with nonnegative coordinates.

\begin{description}
\item[PA1] We have $\gamma:=\sum_{j=1}^q\alpha_j+\sum_{j=1}^q\beta_j<1$.
\item[PA2] The process $\left((Z_t,N^{(t)})\right)_{t\in\Z}$ is stationary, ergodic and adapted to a filtration $\left(\mathcal{F}_t\right)_{t\in\Z}$ such that for all $t\in\Z$, $N^{(t)}$ is independent from $\mathcal{F}_{t-1}$. Moreover, $\E\vert Z_1\vert<\infty$.
\end{description}

\begin{prop}\label{PARX}
\begin{enumerate}
\item
Suppose that Assumptions {\bf PA1-PA2} hold. There then exists a unique non-anticipative, stationary and ergodic process $(Y_t)_{t\in\Z}$ solution of (\ref{Poisss}). 
\item
Additionally, if Assumption {\bf A4} is also satisfied with $\varepsilon_t=N^{(t)}$, there then exists $C>0$ and $\rho\in(0,1)$, such that for all $h\geq 1$, $\Theta_{1,h}\left((Y_t,\lambda_t)_t\right)\leq C \left[\rho^h+\sum_{i=1}^h\rho^i\Theta_{1,h-i}(Z)\right]$. 
\end{enumerate}
\end{prop}

\paragraph{Note.} Our result extends substantially that of \citet{Cav}. First, we prove ergodicity properties in PARX models without assuming that the covariate process $(Z_t)_{t\in\Z}$ is a Markov chain defined by a random map contracting in average. Secondly, for the stochastic dependence properties, we control the coefficient of functional dependence measure only assuming a general Bernoulli shift representation for $(Z_t)_{t\in\Z}$.   
For instance, $(Z_t)_{t\in\Z}$ can be defined by an infinite moving average process and is not necessarily Markovian.

\subsection{Dynamic binary choice model}\label{binary}
We consider the dynamic 
\begin{equation}\label{equa}
Y_t=\mathds{1}_{g(Y_{t-1},\ldots,Y_{t-q},\zeta_t)>0},
\end{equation}
where $(\zeta_t)_{t\in\Z}$ is a stationary process taking values in a measurable space $E'$ and $g:\{0,1\}^q\times E'\rightarrow \R$ is a measurable function. 
This kind of binary model is popular in econometrics for studying the dynamics of recessions. See \citet{deJ} who studied the case $g$ linear and \citet{kauppi} for a study of US recessions.


\begin{prop}\label{supercond}
\begin{enumerate}
\item
Assume that $(\zeta_t)_{t\in \Z}$ is a stationary and ergodic process such that
\begin{equation}\label{good}
\P\left(\min_{y\in\{0,1\}^q, 1\leq t\leq q} g(y,\zeta_t)>0\right)+\P\left(\max_{y\in\{0,1\}^q, 1\leq t\leq q} g(y,\zeta_t)\leq 0\right)>0.
\end{equation}
There then exists a unique stationary and ergodic solution $(Y_t)_{t\in\Z}$ for the recursions (\ref{equa}).
\item
Assume that for some real numbers $a_1,\ldots,a_q$ and $\pi\in\R^e$, $g(y,\zeta_t)=\sum_{i=1}^q a_i y_i+\pi'Z_{t-1}+\varepsilon_t$, with $\zeta_t=(Z_{t-1},\varepsilon_t)$ satisfying {\bf A4} and the c.d.f. $F_{\varepsilon}$ of $\varepsilon_t$ being Lipschitz and taking values in $(0,1)$. 
Moreover, setting $\upsilon_t=\pi'Z_{t-1}+\varepsilon_t$, we assume that there exists $\delta>0$ and a positive integer $K$ such that 
\begin{equation}\label{deJ}
\P\left(\phi_{-}+\min_{1\leq t\leq q}\upsilon_t>0 \vert \mathcal{F}_{-K}\right)+\P\left(\phi_{+}+\max_{1\leq t\leq q}\upsilon_t\leq 0 \vert \mathcal{F}_{-K}\right)\geq \delta \mbox{ a.s.},
\end{equation}
where 
$$\phi_{+}=\max\left\{\sum_{i=1}^q a_iy_i: (y_1,\ldots,y_n)\in\{0,1\}^n\right\},\quad \phi_{-}=\min\left\{\sum_{i=1}^q a_iy_i: (y_1,\ldots,y_n)\in\{0,1\}^n\right\}.$$
There then exists $C>0$ and $\rho\in(0,1)$, such that for all $h\geq 1$, 
$$\Theta_{1,h}(Y)\leq C \left[\rho^h+\sum_{i=1}^h\rho^i\Theta_{1,h-i}(Z)\right].$$
\end{enumerate}
\end{prop}

\paragraph{Notes} 
\begin{enumerate}
\item
Consider the case of $g$ linear as in the second point of Proposition \ref{supercond}. In this case, \citet{deJ} derived existence of a unique stationary and ergodic solution for when Condition (\ref{deJ}) holds true. As shown in \citet{deJ}, Condition (\ref{deJ}) holds in particular when the process $(\upsilon_t)_{t\in\Z}$ is $m-$dependent or for some infinite moving averages. 
Condition (\ref{good}) is much weaker since it holds as soon as 
\begin{equation}\label{oncompare}
\P\left(\phi_{+}+\max_{1\leq t\leq q}\upsilon_t\leq 0\right)+\P\left(\phi_{-}+\min_{1\leq t\leq q}\upsilon_t>0\right)>0.
\end{equation}
Condition (\ref{oncompare}) holds true as soon as the random vector $(\upsilon_1,\ldots,\upsilon_q)$ has full support. 
Another sufficient condition for (\ref{oncompare}) is the following. If $(\upsilon_t)_{t\in\Z}$ is adapted to a filtration $\left(\mathcal{F}_t\right)_{t\in\Z}$, we assume that for any $x\in\R$ and $t\in\Z$, $\P\left(\upsilon_t\leq x\vert\mathcal{F}_{t-1}\right)>0$ a.s. or for any $t\in\Z$ and $x\in \R$, $\P\left(\upsilon_t>x\vert \mathcal{F}_{t-1}\right)>0$ a.s. Recall that $\upsilon_t=\pi'Z_{t-1}+\varepsilon_t$.
The latter condition is valid in particular when $\varepsilon_t$ has full support and is independent from $\mathcal{F}_{t-1}=\sigma\left((\varepsilon_s,Z_s): s\leq t-1\right)$. 

\item 
As the proof of Proposition \ref{supercond} will show, the condition (\ref{deJ}) implies Assumption {\bf A2}. Condition (\ref{good}) is only used for applying Theorem \ref{bibli}. However, (\ref{good}) does not entail mixing properties. In contrast, Condition (\ref{deJ}) does. See \citet{deJ}, Theorem $2$. Our results (see point $2.$ of Proposition \ref{supercond}) give a complement when the covariate process is not necessarily strongly mixing and has a Bernoulli shift representation. 
\item
When $\zeta_t=(Z_{t-1},\varepsilon_t)\in \R^{d+1}$ in (\ref{equa}), one can allow interactions between lag values of the response and the covariates. For example,
$$g(y,\zeta_t)=\sum_{i=1}^d c_i y_i+\sum_{i=1}^q\sum_{j=1}^d \left[a_{i,j}y_i+b_{i,j}(1-y_i)\right] Z_{j,t-i}+\varepsilon_t.$$ 
When $\varepsilon_t$ is independent of $\mathcal{F}_{t-1}=\sigma\left((\varepsilon_{t-j},Z_{t-j}): j\geq 1\right)$, one can show, using the same arguments as in the previous point, that condition (\ref{good}) is satisfied as soon as the distribution of $\varepsilon_t$ has support equal to the whole real line.
 We will not give a control of the functional dependence measure for this model because we were not able to check ${\bf A2}$ when the covariate process $(Z_t)_{t\in\Z}$ is not bounded. However when the cdf of $\varepsilon_t$ is known (e.g. for the logistic or the probit model), it is widely known that ergodicity of the process is sufficient for showing consistency and asymptotic normality of conditional pseudo likelihood estimators of the parameters.
\end{enumerate}

\subsection{Categorical time series with covariates}\label{CATCAT}
We consider a finite set $E=\{1,2,\ldots,N\}$, an integer $q\geq 1$, a process $(Z_t)_{t\in\Z}$ taking values in $\mathcal{Z}\subset\R^d$ and
a family $\left\{K_z\left(\cdot\vert \cdot\right): z\in \mathcal{Z}\right\}$ 
of probability kernels from $E^q$ to $E$. Our aim is to construct a process $(Y_t)_{t\in\Z}$, taking values in $E$ and such that 
$$\P\left(Y_t=i\vert Y_{t-1}^{-},Z_{t-1}^{-}\right)=K_{Z_{t-1}}\left(i\vert Y_{t-1},\ldots,Y_{t-q}\right).$$
A particular example is given by the multinomial autoregression, i.e.
$$K_z(i,y_1,\ldots,y_q)=\frac{\exp\left(\sum_{j=1}^q a_{i,j} y_j+\gamma_i' z\right)}{\sum_{k=1}^N\exp\left(\sum_{j=1}^q a_{k,j} y_j+\gamma_k' z\right)}$$
and is a classical model for categorical time series. See \citet{Fokianos2003}. In econometrics, \citet{russell} studied the dynamic of price changes using this kind of model but with a more general observation-driven form such as in GARCH models and that will not fall into our framework.

For applying our results, we now define some random maps.  
For $t\in\Z$, let $\varepsilon_t$ be a random variable uniformly distributed over $[0,1]$. 
For $u\in [0,1]$, $z\in E_1'$, $y\in E^q$ and $u\in [0,1]$, we set 
$$K_z^{-}(u\vert y)=\inf\left\{i=1,\ldots,N: \sum_{j=1}^i K_z(j\vert y)\geq u\right\}$$
and 
$$f_t(y_1,\ldots,y_q)=\left(K^{-}_{Z_{t-1}}\left(\varepsilon_t\vert y_1,\ldots,y_q\right),y_1,\ldots,y_{q-1}\right)'.$$

We introduce the following assumptions.

\begin{description}
\item[C1]
The probability kernels $K_z$ are lower bounded by a positive constant, i.e. for any $z\in E'$,  
$\eta(z):=\min_{(i,y)\in E^{q+1}}K_z\left(i\vert y\right)>0$. 
\item[C2]
The process $\left((Z_t,\varepsilon_t)\right)_{t\in\Z}$ is stationary and ergodic. Moreover,
for $t\in\Z$, $\varepsilon_t$ is independent from $\mathcal{F}_{t-1}=\sigma\left((Z_s,\varepsilon_s): s\leq t-1\right)$.

\item[C3] There exists a constant $C>0$ such that for all $y_1,\ldots,y_q\in E$, 
$$\sum_{i=1}^N\left\vert K_{z}(i\vert y_1,\ldots,y_q)-K_{\overline{z}}(i\vert y_1,\ldots,y_q)\right\vert \leq C\vert z-\overline{z}\vert.$$
\end{description}

\begin{prop}\label{cat}
Suppose that Assumptions {\bf C1-C2} hold.
\begin{enumerate}
\item
 There exists a unique stationary process satisfying the recursions 
\begin{equation}\label{recurcur}
Y_t=K_{Z_{t-1}}^{-}\left(\varepsilon_t\vert Y_{t-1},\ldots,Y_{t-q}\right),\quad t\in\Z.
\end{equation}
Moreover, the process $\left((Y_t,Z_t)\right)_{t\in\Z}$ is ergodic.
\item
Additionally, assume that Assumption {\bf A4} and {\bf C3} hold true and that $\eta_{-}=\inf_{z\in E_1'}\eta(z)>0$ in {\bf C1}. 
There then exist $C>0$ and $\rho\in(0,1)$, such that for all $h\geq 1$, $\Theta_{1,h}(Y)\leq C \left[\rho^h+\sum_{i=1}^h\rho^i\Theta_{1,h-i}(Z)\right]$.
\end{enumerate}
\end{prop}

\paragraph{Note.} A proof of the first point of Proposition \ref{cat} is based on Theorem \ref{bibli} and provides a general result for existence of stationary categorical time series with covariates. In particular, probit, logistic and multinomial autoregressions can be considered without restriction for the covariate process $(Z_t)_{t\in\Z}$. However, the derivation of the dependence properties in the second point imposes a more restrictive assumption on the transition kernel $K$ because it is necessary to check Assumption {\bf A2}.   
Recently, \citet{FT1} studied categorical time series under the strict exogeneity assumption for the covariate process. 
For the recursions (\ref{recurcur}), strict exogeneity holds true as soon as the two processes $(Z_t)_{t\in\Z}$ and $(\varepsilon_t)_{t\in\Z}$ are independent. 
Assumption {\bf C2} is weaker than strict exogeneity in general.

\subsection{Categorical time series and coalescence of the paths}
\begin{figure}[H]
\begin{center}
\definecolor{ffqqqq}{rgb}{1,0,0}\definecolor{ududff}{rgb}{0.30196078431372547,0.30196078431372547,1}\definecolor{xdxdff}{rgb}{0.49019607843137253,0.49019607843137253,1}
\begin{tikzpicture}[line cap=round,line join=round,>=triangle 45,x=0.78cm,y=0.78cm]\clip(-15.5,-5) rectangle (6,4);

\draw [green,line width=1pt] (-5,-0.1) -- (-1.5,2);
\draw [green,line width=1pt] (-12,1.95) -- (-8.5,1.95);
\draw[green,line width=1pt] (-12,-0.1)--(-8.5,2);
\draw [green,line width=1pt] (-8.5,2) -- (-5,-0.1);
\draw [green,line width=1pt] (-12,-2.1) -- (-8.5,2);

\draw [line width=0.5pt] (-12,2) -- (-8.5,2);
\draw [line width=0.5pt] (-8.5,2) -- (-5,2);
\draw [line width=0.5pt] (-5,2) -- (-1.5,2);
\draw [line width=0.5pt] (-1.5,2) -- (-0.5,2);
\begin{scriptsize}
\draw [fill=xdxdff] (-12,2) circle (2.5pt);
\draw [fill=xdxdff] (-8.5,2) circle (2.5pt);
\draw [fill=xdxdff] (-5,2) circle (2.5pt);
\draw [fill=xdxdff] (-1.5,2) circle (2.5pt);
\end{scriptsize}

\draw [line width=0.5pt] (-12,0) -- (-8.5,0);
\draw [line width=0.5pt] (-8.5,0) -- (-5,0);
\draw [line width=0.5pt] (-5,0) -- (-1.5,0);
\draw [line width=0.5pt] (-1.5,0) -- (-0.5,0);
\begin{scriptsize}
\draw [fill=xdxdff] (-12,0) circle (2.5pt);
\draw [fill=xdxdff] (-8.5,0) circle (2.5pt);
\draw [fill=xdxdff] (-5,0) circle (2.5pt);
\draw [fill=xdxdff] (-1.5,0) circle (2.5pt);
\end{scriptsize}

\draw [line width=0.5pt] (-12,-2) -- (-8.5,-2);
\draw [line width=0.5pt] (-8.5,-2) -- (-5,-2);
\draw [line width=0.5pt] (-5,-2) -- (-1.5,-2);
\draw [line width=0.5pt] (-1.5,-2) -- (-0.5,-2);
\begin{scriptsize}
\draw [fill=xdxdff] (-12,-2) circle (2.5pt);
\draw [fill=xdxdff] (-8.5,-2) circle (2.5pt);
\draw [fill=xdxdff] (-5,-2) circle (2.5pt);
\draw [fill=xdxdff] (-1.5,-2) circle (2.5pt);
\end{scriptsize}

\begin{scriptsize}
\draw[color=red](-14,1.9) node {State $1$};
\draw[color=red](-14,-0.1) node {State $2$};
\draw[color=red](-14,-2) node {State $3$};
\draw[color=blue](-14,-2.75) node {Time };
\draw[color=blue](-12,-2.75) node {$t-3$};
\draw[color=blue](-8.5,-2.75) node {$t-2$};
\draw[color=blue](-5,-2.75) node {$t-1$};
\draw[color=blue](-1.5,-2.75) node {$t$};
\end{scriptsize}

\end{tikzpicture}
\end{center}
\caption{Illustration of the convergence for $N=3$ modalities and $q=1$ lag \label{fig1}}
\end{figure}
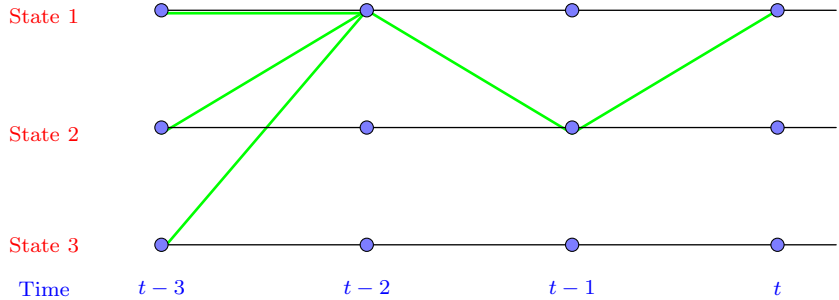

In this section, we give another interpretation of the convergence of the backward iterations for categorical time series.
This interpretation has a link with some perfect simulation schemes that are widely known for Markov chains. See \citet{PW}.
Since the state space is discrete, the iterations should be automatically constant after some steps.
Figure \ref{fig1} illustrates the convergence when $q=1$ and $N=3$. In this case, $f_t(j)=K^{-}_{Z_{t-1}}(\varepsilon_t\vert j)$ for $j=1,2,3$.
Setting 
$$T=\inf\left\{k\geq 1: \varepsilon_{t-k}\leq \min_{i,j}K_{Z_{t-k-1}}(i\vert j)\right\},$$
we know that from the ergodicity assumption {\bf C2} and the positivity assumption {\bf C1}, $T$ is finite almost surely.
In this case, $f_{t-T}(j)=1$ a.s., $f_{-\infty}^t:=\lim_{n\rightarrow \infty}f_{t-n}^t(j)=f_{t-T}^t(j)$ a.s. and the limit does not depend on the state $j$.
All the paths corresponding to $f_{t-n}^t(j)$ for $n\geq T$ coalesce through State $1$.

When $q\geq 2$, our assumptions guarantee that it is possible to get $q$ times successively the value $1$ for the time series whatever the previous values.
A coalescence property for the paths will then also occur in this case.

This interpretation is also relevant for getting an interpretation of Condition (\ref{good}) for binary choice models.
When (\ref{good}) is satisfied, it is possible to get, starting at any time $t$, either $q$ times the value $1$ or $q-$times the value $0$ whatever the previous values of the binary time series. Running the backward iterations, we have coalescence of the paths at the first (random) time $t-T$ such that such an event occurs.

\section{Conclusion}\label{Conc}
A general theoretical analysis of nonlinear autoregressive time series models with exogenous covariates is absent from the present time series literature and only a few references consider such a problem, mainly for specific examples. The aim of this paper was to provide some results for a reasonable class of nonlinear time series models for which the required assumptions can be checked. 
In particular, we provide two results, Theorem \ref{main} and Theorem \ref{autoreg}, which justify existence and uniqueness 
of stationary and ergodic solutions possessing some moments. The crucial assumption to check, {\bf A2} or {\bf B3}, involves a uniform conditional contraction condition. Assumption {\bf B3} is the main assumption to check for autoregressive models with several lags.
For some nonlinear models already considered in the literature, such as GARCH or autoregressive Poisson processes, this contraction condition is easily checked because the exogenous covariates have an additive contribution in the expression of the latent process and play the role of a random intercept which does not modify the usual stability conditions. 
However, our results can also be applied to autoregressive processes for which lag parameters depend on the covariates (see Section \ref{comp}, the notes after Theorem \ref{autoreg}, Proposition \ref{CHARME} and Proposition \ref{supercond} for some examples). 
In this case, a uniform control of the random lag parameters is necessary to check our assumptions which shows the limit of our approach.
 
It may be possible to weaken our uniform contraction condition, as shown in Theorem \ref{relax}, at least under a strict exogeneity assumption. However, getting additional general results to ensure existence of some unconditional moments, to control dependence coefficients and to consider higher-order autoregressive processes would require substantial effort.
The proposed framework is also useful for deriving weak dependence properties of the solution, leading to the possibility to apply many existing statistical inference procedures, the central limit theorem of Section \ref{TLC} providing an illustration. To this end, the functional dependence measure discussed in Section \ref{mix} is of primary importance.
Note that a general result for getting weak dependence properties of autoregressive processes with exogenous covariates is 
also new and it is another contribution of this paper. Finally, we also derived results for categorical time series in Sections \ref{binary}-\ref{CATCAT}. Apart from the weak dependence properties which can be derived from our general results, 
we also obtained stationarity conditions with weaker assumptions, applying Theorem \ref{bibli}. Note that whatever the results used in the paper (Theorems \ref{bibli}, \ref{main} or \ref{relax}), the convergence of the backward iterations of random maps appears to be a central point of view for considering many different types of autoregressive systems with exogenous regressors.

\section{Proofs of the results}
\subsection{Proof of Theorem \ref{main}}
We use the convention $f_t^{t-1}(x)=x$ for $(x,t)\in E\times \Z$.
From Assumption {\bf A2}, if $(t,s,s')\in \Z^3$ are such that $s'\leq s\leq t$ and $x,y\in E$,  
$$\E\left[d^p\left(f_s^t(x),f_{s'}^t(y)\right)\vert \mathcal{F}_{t-1}\right]\leq L^p d^p\left(f_s^{t-1}(x),f_{s'}^{t-1}(y)\right),$$
and then  
\begin{equation}\label{bound1}
\Vert d\left(f_s^t(x),f_{s'}^t(y)\right)\Vert_p\leq L\Vert d\left(f_s^{t-1}(x),f_{s'}^{t-1}(y)\right)\Vert_p.
\end{equation}
Applying (\ref{bound1}) with $y=x_0$, $s=s'=t$, we get $\sup_{t\in\Z}\Vert d\left(f_t(x),f_t(x_0)\right)\Vert_p<\infty$.
Next, using {\bf A1} and the triangular inequality, we get 
\begin{equation}\label{bound2}
\sup_{t\in \Z}\Vert d\left(f_t(x),y\right)\Vert_p<\infty\mbox{ for every }(x,y)\in E^2.
\end{equation}

With the same kind of arguments, we get for $s\leq t-m$,  
\begin{equation}\label{bound22}
\Vert d\left(f_s^t(x),f_{s'}^t(y)\right)\Vert_p\leq \kappa\Vert d\left(f_s^{t-m}(x),f_{s'}^{t-m}(y)\right)\Vert_p.
\end{equation}
\begin{enumerate}
\item
We denote by $[z]$ the integer part of a real number $z$. To apply recursively the previous bounds, we note that for any integer $i$,
$i+1=r_1m+r_2$ with $r_1=[(i+1)/m]$ and $r_2=i+1-r_1m$.
We then get from (\ref{bound1}) and (\ref{bound22}), 
\begin{eqnarray*}
\sum_{i\geq 0}\Vert d\left(f_{t-i}^t(x),f_{t-i-1}^t(x)\right)\Vert_p&\leq& 
\sum_{i\geq 0}\kappa^{[(i+1)/m]}L^{i+1-m[(i+1)/m]}\Vert d\left(x,f_{t-i-1}(x)\right)\Vert_p \\
&\leq& \frac{(L+1)^m\kappa^{(1-m)/m}}{1-\kappa^{1/m}}\sup_{j\in\Z}\Vert d\left(x,f_j(x)\right)\Vert_p.
\end{eqnarray*}
This latter bound entails that the series $\sum_{i\geq 0}d\left(f_{t-i}^t(x),f_{t-i-1}^t(x)\right)$ is almost surely finite. By the Cauchy criterion, there exists a random variable 
$X_t(x)$ such that $\lim_{i\rightarrow \infty} d\left(f_{t-i}^t(x),X_t(x)\right)=0$ a.s.
Moreover, from the previous bound, (\ref{bound2}) and the triangular inequality, we deduce that 
\begin{equation}\label{bound3}
\sum_{s,t\in\Z, s\leq t}\Vert d\left(y,f_s^t(x)\right)\Vert_p<\infty\mbox{ for every }(x,y)\in E^2.
\end{equation}
Next we note that the convergence also holds in $\L^p$, since from Fatou's lemma,
\begin{eqnarray*}
\Vert d\left(X_t(x),f_{t-s}^t(x)\right)\Vert_p&\leq& \liminf_{j\rightarrow\infty}\Vert d\left(f_{t-j}^t(x),f_{t-s}^t(x)\right)\Vert_p\\
&\leq& \sum_{i\geq s}\Vert d\left(f_{t-i}^t(x),f_{t-i-1}^t(x)\right)\Vert_p \\
 & \leq &   \kappa^{s/m}\frac{(L+1)^m\kappa^{(1-m)/m}}{1-\kappa^{1/m}}\sup_{j\in\Z}\Vert d\left(f_j(x),x\right)\Vert_p \rightarrow 0 \text{ as } s \rightarrow \infty.
\end{eqnarray*}
Finally, using (\ref{bound3}) and the triangular inequality, we get the last assertion $\sup_{t\in\Z}\Vert d\left(X_t(x),x_0\right)\Vert_p<\infty$.  

\item
If $x\neq y$, we have from the almost sure convergence and Fatou's lemma
\begin{eqnarray*}
\Vert d\left(X_t(x),X_t(y)\right)\Vert_p&\leq& \liminf_{s\rightarrow \infty}\Vert d\left(f_{t-s}^t(x),f_{t-s}^t(y)\right)\Vert_p\\
&\leq&\liminf_{s\rightarrow \infty}\kappa^{[(s+1)/m]}L^{s+1-[(s+1)/m]m}d(x,y)=0.
\end{eqnarray*}
This shows the second point.

\item
For the third point, we observe that for any $j\geq 1$, there exists a measurable function $H_j^{(x)}:E'^{j+1}\rightarrow E$ such that 
$f_{t-j}^t(x)=H_j^{(x)}\left(\zeta_t,\ldots,\zeta_{t-j}\right)$. Since $\lim_{j\rightarrow\infty}H_j^{(x)}$ exists $\P_{(\zeta_{t-j})_{j\geq 0}}$ a.s.,
it is then possible to define a measurable function $H:(E')^{\N}\rightarrow E$ such that 
$X_t=H\left((\zeta_{t-j})_{j\geq 0}\right)$ a.s.
The process $(X_t)_{t\in\Z}$ has a Bernoulli shift structure with dependent entries and is then stationary and ergodic provided that the process $(\zeta_t)_{t\in\Z}$ satisfies the same properties. 

\item
The last property follows from the following bounds which hold for any $j\geq 1$:
\begin{eqnarray*}
\Vert d(X_t,Y_t)\Vert_p&=& \Vert d\left(f_{t-mj+1}^t(X_{t-mj}),f_{t-mj+1}^t(Y_{t-mj}\right)\Vert_p\\
&\leq& \kappa^j\left[\sup_{j\in\Z}\Vert d(x_0,X_j)\Vert_p+\sup_{j\in\Z}\Vert d(x_0,Y_j)\Vert_p\right].  \quad   \quad \square
\end{eqnarray*}
\end{enumerate}

\subsection{Proof of Theorem \ref{relax}}

\begin{enumerate}
\item
From Assumption {\bf A0}, $f_{t-s}^t(x)$ is, conditionally on $Z$, an iteration of $s+1$ independent random maps.
Using Assumption {\bf A2'}, we get
$$\sum_{s\geq 0}\E^{1/p}\left[d^p\left(f_{t-s}^t(x),f_{t-s-1}^t(x)\right)\vert Z\right]\leq \sum_{s\geq 0}\prod_{i=1}^{s+1}\kappa(Z_{t-i})b_{t-s-2}(x),$$
with 
$$b^p_t(x)=\int d^p\left(x,F(x,Z_t,u)\right)d\P_{\varepsilon_0}(u).$$
From Assumption {\bf A2'} and the triangular inequality, we have $\E\log^{+}b_t(x)<\infty$ for any $x\in E$.
We are going to show that
\begin{equation}\label{ncv}
\sum_{s\geq 0}\E^{1/p}\left[d^p\left(f_{t-s}^t(x),f_{t-s-1}^t(x)\right)\vert Z\right]<\infty\mbox{ a.s.}
\end{equation}
This follows from the assumptions on the logarithmic moments. Indeed, $\left((\kappa(Z_t),b_t(x))\right)_{t\in\Z}$ is a stationary process and it is widely known that the stochastic recursions
$$Y_t=\kappa(Z_{t-1})Y_{t-1}+b_{t-1}(x)$$
have a unique stationary solution given by 
$$Y_t=b_{t-1}(x)+\sum_{s\geq 0}\prod_{i=1}^{s+1}\kappa(Z_{t-i})b_{t-s-2}(x),$$
the latter series being convergent almost surely. See for instance \citet{Brandt}, Theorem $1$. This shows (\ref{ncv}).
Using Minkowski's inequality and Fatou's lemma for conditional expectations, see for instance \citet{Kall}, Chapter $5$, we then deduce that
$$\E^{1/p}\left[S(x)^p\vert Z\right]<\infty\mbox{ a.s. } S(x):=\sum_{s\geq 0}d\left(f_{t-s}^t(x),f_{t-s-1}^t(x)\right).$$
As a consequence, we have $\P\left(S(x)<\infty\vert Z\right)=1$ a.s. and then 
$\P\left(S(x)<\infty\right)=1$. From the Cauchy criterion, we then conclude the existence of a random variable $X_t(x)$ such that 
$\lim_{s\rightarrow \infty}f_{t-s}^t(x)=X_t(x)$ a.s. 
Note that from Assumption {\bf A1'} and (\ref{ncv}), we have $\E\left[d^p\left(f_{t-s}^t(x),x_0\right)\vert Z\right]<\infty$ a.s. for every positive integer $s$.
The convergence $\lim_{s\rightarrow \infty}\E\left[d^p\left(f_{t-s}^t(x),X_t(x)\right)\vert Z\right]=0$ a.s. and $\E\left[d^p\left(X_t(x),x_0\right)\vert Z\right]<\infty$ a.s.
follow as in the proof of Theorem \ref{main}, using Fatou's lemma for conditional expectation.

\item
For a positive integer $s$, we have from {\bf A2}, 
$$\E^{1/p}\left[d^p\left(f_{t-s}^t(x),f_{t-s}^t(y)\right)\vert Z\right]\leq \prod_{i=1}^s\kappa(Z_{t-i})d(x,y)\rightarrow 0\mbox{ a.s}.$$
Letting $s\rightarrow \infty$, we have $\E\left[d^p(X_t(x),X_t(y))\vert Z\right]=0$ a.s. and then $\P\left(X_t(x)\neq X_t(y)\vert Z\right)=0$ a.s. Taking the expectation, we conclude that $\P\left(X_t(x)\neq X_t(y)\right)=0$.

\item
From the almost sure convergence of the sequence $\left(f_{t-s}^t(x)\right)_{s\geq 0}$, stationarity and ergodicity of the process $\left((X_t,Z_t)\right)_{t\in\Z}$ follows exactly as in the proof of point $3$ of Theorem \ref{main}.

\item
Let $(Y_t)_{t\in\Z}$ be a stochastic process satisfying the proposed conditions. 
If the process is non-anticipative, we have from {\bf A2'}, 
\begin{eqnarray*}
\E\left[d^p(X_t,Y_t)\vert Z\right]&=&\E\left[d^p\left(f_{t-s}^t(X_{t-s-1}),f_{t-s}^t(Y_{t-s-1})\right)\vert Z\right]\\
&\leq& \prod_{i=0}^s\kappa^p(Z_{t-i-1})\cdot\E\left[d^p\left(X_{t-s-1},Y_{t-s-1}\right)\vert Z\right].
\end{eqnarray*}
Note that, from {\bf A2'}, $\prod_{i=0}^s\kappa^p(Z_{t-i-1})=o_{\P}(1)$. Moreover 
$$\E^{1/p}\left[d^p\left(X_{t-s-1},Y_{t-s-1}\right)\vert Z\right]\leq \E^{1/p}\left[d^p\left(X_{t-s-1},x_0\right)\vert Z\right]+\E^{1/p}\left[d^p\left(x_0,Y_{t-s-1}\right)\vert Z\right].$$
Note that 
$$\E^{1/p}\left[d^p\left(x_0,Y_{t-s-1}\right)\vert Z\right]=\E^{1/p}\left[d^p\left(x_0,Y_{t-s-1}\right)\vert Z_{t-s-1},Z_{t-s-2},\ldots\right]$$
and if the process $\left((Z_t,Y_t)\right)_{t\in\Z}$ is stationary, then the process $(V_t)_{t\in\Z}$ defined by $$V_t=\E^{1/p}\left[d^p\left(x_0,Y_t\right)\vert Z_t,Z_{t-1},\ldots\right]$$
 is also stationary and takes finite values from our assumptions. Then $V_{t-s-1}=O_{\P}(1)$ where $s$ is the index of the sequence. The same property holds if $Y_t$ is replaced with $X_t$. 
As a consequence 
$$\E\left[d^p\left(X_{t-s-1},Y_{t-s-1}\right)\vert Z\right]=O_{\P}(1).$$
We then conclude that $\E\left[d^p(X_t,Y_t)\vert Z\right]=0$ a.s. Then $\P\left[X_t\neq Y_t\vert Z\right]=0$ a.s. and by integration, 
we get the conclusion.

\end{enumerate}

\subsection{Proof of Proposition \ref{WD}}
\begin{enumerate}
\item
We use the decomposition
\begin{eqnarray*} 
\overline{f}_{-\infty}^t-f_{-\infty}^t&=&
\sum_{i=0}^{t-1}\left[f_{t-i}^t\circ \overline{f}_{-\infty}^{t-i-1}-f_{t-i-1}^t\circ \overline{f}_{-\infty}^{t-i-2}\right]\\
&+& \overline{f}_t\circ \overline{f}_{-\infty}^{t-1}-f_t\circ \overline{f}_{-\infty}^{t-1}.
\end{eqnarray*}
From Assumption {\bf A2} and Assumption {\bf A4}, we have, for $i=0,\ldots,t-2$,
\begin{eqnarray*}
\Vert d\left(f_{t-i}^t\circ \overline{f}_{-\infty}^{t-i-1},f_{t-i-1}^t\circ \overline{f}_{-\infty}^{t-i-2}\right)\Vert_p&\leq& \kappa^{\frac{i+1}{m}-1} L^m\Vert d\left(\overline{f}_{t-i-1}\circ \overline{f}^{t-i-2}_{-\infty},f_{t-i-1}\circ \overline{f}^{t-i-2}_{-\infty}\right)\Vert_p\\
&\leq& \kappa^{\frac{i+1}{m}-1} L^m \Vert S(\overline{f}_{-\infty}^{t-i-2}) H_{t-i-2}\Vert_p\\
&\leq& \kappa^{\frac{i+1}{m}-1} L^m \Vert S(f_{-\infty}^0)\Vert_s \Vert H_{t-i-2}\Vert_r.
\end{eqnarray*}
If $i=t-1$, we have $\Vert d\left(f_{t-i}^t\circ \overline{f}_{-\infty}^{t-i-1},f_{t-i-1}^t\circ \overline{f}_{-\infty}^{t-i-2}\right)\Vert_p\leq 2\Vert d\left(0,f_{-\infty}^0\right)\Vert_p\kappa^{t/m-1}$. 
Using the triangular inequality, we get for $t\geq 2$,
$$\theta_{p,t}\leq \kappa^{-1} L^m \Vert S(f_{-\infty}^0)\Vert_s\sum_{i=0}^{t-2}\kappa^{(i+1)/m}\Vert H_{t-i-2}\Vert_r+\Vert S(f_{-\infty}^0)\Vert_s\Vert H_{t-1}\Vert_r+2\kappa^{t/m-1}\Vert d\left(0,f_{-\infty}^0\right)\Vert_p.$$
The bound (\ref{bbound}) is obtained by summation and entails the simpler bound (\ref{goodmaj}).

\item
From {\bf A4'}, we have {\bf A4} with $H_{t-1}=d\left(Z_{t-1},\overline{Z}_{t-1}\right)$ with 
$$\overline{Z}_t=H'\left(\eta_t,\ldots,\eta_1,\eta_0',\eta_{-1},\ldots\right).$$
We then deduce the result from the previous point, noticing that $\eta_{r,j}=\Theta_{r,j-1}(Z)$.
\end{enumerate}

\subsection{Proof of Theorem \ref{autoreg}}
Define the following random map 
$$f_t(u_1,\ldots,u_q)=\left(F\left(u_1,\ldots,u_q,\zeta_t\right)',u'_1,\ldots,u'_{q-1}\right)'.$$
We set $x=(u_1,\ldots,u_q)\in E^q$ and for $1\leq t\leq q$, $U_t(x)=u_{q-t+1}$. Next for $t\geq q+1$, we define $U_t(x)$ recursively by 
$$U_t(x)=F\left(U_{t-1}(x),\ldots,U_{t-q}(x),\varepsilon_t\right).$$ 
We then have for $t\geq q+1$,
$$\left(U_t(x),\ldots,U_{t-q+1}(x)\right)=f_{q+1}^t(x).$$
Using our assumptions, we have for $t\geq q+1$, 
$$\E\left[\Delta^p_{vec}\left(U_t(x),U_t(x')\right)\vert \mathcal{F}_{t-1}\right]\preceq \sum_{i=1}^q A_i \Delta^p_{vec}\left(U_{t-i}(x),U_{t-i}(x')\right).$$
We introduce the matrix
$$B=\begin{pmatrix} A_1&\cdots& A_{q-1}& A_q\\
& I_{k(q-1)}& & 0_{k(q-1),1} \end{pmatrix}.$$
The condition $\rho(A_1+\cdots+A_q)<1$ entails that $\rho(B)<1$. Indeed, if $v=(v_1',\ldots,v'_q)'\in \R^{kq}\setminus\{0\}$ is such that $B v=\lambda v$ for $\vert\lambda\vert\geq 1$, we get the equality 
$v_1=\left[\lambda^{-1}A_1+\cdots+\lambda^{-q}A_q\right]v_1$. Since the coefficients of the $A_i'$s are nonnegative, we get 
$$\vert v_1\vert_{vec}\preceq \sum_{i=1}^q\vert\lambda\vert^{-i}A_i\vert v_1\vert_{vec}\preceq \sum_{i=1}^q A_i\vert v_1\vert_{vec},$$
where $\vert v_1\vert_{vec}$ denotes the vector of the absolute values of the coordinates of $v_1$.
We then get $\vert v_1\vert_{vec}\preceq \left(\sum_{j=1}^q A_j\right)^k\vert v_1\vert_{vec}$ for any positive integer $k$. Letting $k\rightarrow \infty$, we obtain $v_1=0$. Since $v_i=\lambda v_{i+1}$ for $i=1,\ldots,q-1$, we get $v=0$ which is a contradiction. Then $\vert \lambda\vert <1$ and $\rho(B)<1$.
Next, we set 
$$V_t(x)=\left(U_t(x)',\ldots,U_{t-q+1}(x)'\right),\quad t\geq q+1.$$
Note that $V_t(x)=f_{q+1}^t(x)$.
We then have 
$$\E\left[\Delta^p_{vec}\left(V_t(x),V_t(x')\right)\vert \mathcal{F}_{t-1}\right]\preceq B \Delta^p_{vec}\left(V_{t-1}(x),V_{t-1}(x')\right)\preceq\cdots\preceq B^{t-q}\Delta^p\left(x,x'\right).$$
Setting for $v,v'\in \R^{kq}$, $\overline{d}(v,v')=\left(\sum_{i=1}^{kq}\Delta^p(v_i,v_i')\right)^{1/p}$,
we get 
$$\E\left[\overline{d}^p\left(f_{q+1}^t(x),f_{q+1}^t(x')\right)\vert \mathcal{F}_q\right]\leq \left\vert\mathds{1}'B^{t-q}\right\vert_{\infty} \overline{d}^p(x,x'),$$
where $\mathds{1}$ denotes the vector of $\R^{kq}$ having all its components equal to $1$, and $\vert\cdot\vert_{\infty}$ the infinite norm in $\R^{kq}$.
Since $\rho(B)<1$, if $t$ is large enough, we have $\left\vert \mathds{1}'B^{t-q}\right\vert_{\infty}<1$.
We then conclude that {\bf A2} is satisfied with $d=\overline{d}$, $m=\inf\{j\geq 1:\left\vert \mathds{1}'B^j\right\vert_{\infty}<1\}$ and $\kappa^p=\left\vert \mathds{1}'B^m\right\vert_{\infty}$. 
Moreover, {\bf A1} is a direct consequence of {\bf B2}. The result then follows from Theorem \ref{main}.$\square$

\subsection{Proof of Proposition \ref{FDM}}
Defining 
$$f_t\left(x_1,\ldots,x_q\right)=\left(F\left(x_1,\ldots,x_q,Z_{t-1},\varepsilon_t\right),x_1,\ldots,x_{q-1}\right),$$
we set for $t\geq 0$, $\overline{Z}_t=H\left(\eta_t,\ldots,\eta_1,\eta_0',\eta_{-1},\ldots\right)$.
Using Assumptions {\bf A5-B5}, we have, for $t\geq 1$, 
$$\E^{1/p}\left[d^p\left(f_t(x_1,\ldots,x_q),\overline{f}_t(x_1,\ldots,x_q)\right)\vert \mathcal{F}_{t-1}\vee\sigma(\xi_0')\right]\leq S(x_1,\ldots,x_q)d\left(Z_{t-1},\overline{Z}_{t-1}\right).$$

We then check Assumption {\bf A4}, setting $H_{t-1}=d\left(Z_{t-1},\overline{Z}_{t-1}\right)$.
One can then apply Proposition \ref{WD} using for instance the $\ell_1-$distance on $E^q$ defined by $\overline{d}(u,v)=\sum_{i=1}^q d(u_i,v_i)$ which is equivalent to the distance $\overline{d}$ used in the proof of Theorem \ref{autoreg}. 
Since $\E^{1/r}\left[d^r\left(Z_{t-1},\overline{Z}_{t-1}\right)\right]=\theta_{r,{t-1}}(Z)$, the proof of the proposition is now complete.

\subsection{Proof of Theorem \ref{TLC}}
The proof is a consequence of Theorem $3$ in \citet{Wu}. It is simply necessary to prove that 
$$\Gamma(Y)=\sum_{t=0}^{\infty}\Vert Y_t-\overline{Y}_t\Vert_2<\infty.$$
To this end, for $t\geq 0$, let $\overline{X}_t$ be the random variables obtained by replacing $\xi_0$ by an independent copy $\xi_0'$ in its Bernoulli shift representation. We then have $\overline{Y}_t=f\left(\overline{X}_t,\ldots,\overline{X}_{t-k}\right)$
for $t\geq k$. From the assumption on the function $f$, stationarity and H\"older's inequality, we have, setting $\ell_1=2p/(p-2)$, 
$$\theta_{2,t}(Y)\leq C\left(1+2k\Vert \vert X_0\vert^{o\ell}\Vert_{\ell_1}\right)\cdot\sum_{i=0}^k\theta_{p,t-i}(X).$$
The required condition easily follows by summation.

\subsection{Proof of Proposition \ref{CHARME}}
\begin{enumerate}
\item
Let $d$ be the distance induced by the $\ell_1-$norm on $\R^q$, i.e. $d(x,y)=\sum_{i=1}^q \vert x_i-y_i\vert$.
We use the notation $\vert x-y\vert$ instead of $d(x,y)$. For a square matrix $A$ of size $q\times q$, 
we denote by $\Vert A\Vert$ the corresponding operator norm of $A$. 
We define the sequence of random maps as follows:
$$g_t(x)=\left(f_1(x_1,\ldots,x_q,Z_{t-1})+\varepsilon_t f_2(x_1,\ldots,x_q,Z_{t-1}),x_1,\ldots,x_{q-1}\right)'.$$
We then have 
$$\left\vert g_t(x)-g_t(y)\right\vert_{vec}\preceq \left(\sum_{i=1}^q c_{i,t}\vert x_i-y_i\vert,\vert x_1-y_1\vert,\ldots,\vert x_{q-1}-y_{q-1}\vert\right)'=A_t\cdot\vert x-y\vert_{vec}.$$
Iterating the previous bound, we get for any positive integer $t$,
$$\left\vert g_1^t(x)-g_1^t(y)\right\vert_{vec}\preceq A_t\cdots A_1\cdot\vert x-y\vert_{vec}.$$
We then deduce that $c(g_1^t)\leq \Vert A_t\cdots A_1\Vert$. The result is then a consequence of Theorem \ref{bibli}, using the condition $\chi({\bf A})<0$.

\item
We check the assumptions of Theorem \ref{autoreg}.
First note that from Assumption {\bf CH1}, the process $\zeta_t=(Z_{t-1},\varepsilon_t)$ is ergodic. This entails {\bf B1}. Next we set 
$$F(x_1,\ldots,x_q,\zeta_t)=f_1(x_1,\ldots,x_q,Z_{t-1})+\varepsilon_t f_2(x_1,\ldots,x_q,Z_{t-1}).$$
Our assumptions guarantee that $F\left(\overline{x},\zeta_1\right)\in\L^p$ and, using {\bf CH2}, 
we deduce that $F(x,\zeta_1)\in\L^p$ for any $x\in\R^q$.
This shows {\bf B2}.
Finally, we check {\bf B3}. To this end, for $i=1,\ldots,q$, we set $\delta_i=\sup_z\Vert a_{i,1}(z)+a_{i,2}(z)\vert\varepsilon_1\vert\Vert_p$.

Using Minkowski's inequality for conditional expectations (see for instance \citet{Doob}, Chapter XI, Section $3$), we have 
\begin{eqnarray*}
\E^{1/p}\left[\left\vert F(x_1,\ldots,x_q,\zeta_t)-F(y_1,\ldots,y_q,\zeta_t)\right\vert^p\vert \mathcal{F}_{t-1}\right]&\leq& \E^{1/p}\left[\left(\sum_{i=1}^q c_{i,t}\vert x_i-y_i\vert\right)^p\vert \mathcal{F}_{t-1}\right]\\
&\leq& \sum_{i=1}^q \E^{1/p}\left[c_{i,t}^p\vert \mathcal{F}_{t-1}\right]\cdot\vert x_i-y_i\vert\\
&\leq& \sum_{i=1}^q \delta_i\vert x_i-y_i\vert.  
\end{eqnarray*}
Next using convexity, we get
$$\E\left[\left\vert F(x_1,\ldots,x_q,\zeta_t)-F(y_1,\ldots,y_q,\zeta_t)\right\vert^p\vert \mathcal{F}_{t-1}\right]\leq \left(\sum_{i=1}^q \delta_i\right)^{p-1}\sum_{i=1}^q \delta_i\vert x_i-y_i\vert.$$
{\bf B3} is then a consequence of {\bf CH3}. 

\item
We apply Proposition \ref{FDM}. From the previous points, it is only required to check {\bf B4} which is a consequence of Assumption {\bf CH3}.
\end{enumerate}

\subsection{Proof of Proposition \ref{Garch}}
For the first part, we apply Theorem \ref{autoreg}. To this end, 
we set $E=\R_+^3$, $F=(F_1,F_2,F_3)$, $F_2(y_1,\ldots,y_q,\zeta_t)=(\varepsilon^+_t)^{\delta}y_{1,1}$, 
$F_3(y_1,\ldots,y_q,\zeta_t)=(\varepsilon^{-}_t)^{\delta}y_{1,1}$ and 
$$F_1\left(y_1,\ldots,y_q,\zeta_t\right)=\pi'Z_{t-1}+\sum_{j=1}^q\beta_j y_{1,j}+\left(\alpha_{1+}(\varepsilon_t^+)^{\delta}+\alpha_{1-}(\varepsilon_t^{-})^{\delta}\right)y_{1,1}+\sum_{j=2}^q \left(\alpha_{j+}y_{2,j}+\alpha_{j-}y_{3,j}\right).$$
 We then deduce that Assumption {\bf B3} holds true with 
$$A_1=\begin{pmatrix} \beta_1+\alpha_{1+}s_{+}+\alpha_{1-}s_{-}& 0&0\\
s_+&0&0\\s_{-}&0&0\end{pmatrix},\quad A_j=\begin{pmatrix} \beta_j& \alpha_{j+}&\alpha_{j-}\\
0&0&0\\0&0&0\end{pmatrix},\quad j\geq 2.$$
It is straightforward to show that the matrix $\Gamma:=\sum_{j=1}^q A_j$ has eigenvalues $0$ and $\frac{a\pm \sqrt{a^2+4(bd+ce)}}{2}$ with $a=\sum_{j=1}^q\beta_j+\alpha_{1+}s_{+}+\alpha_{1-}s_{-}$, $b=\sum_{j=2}^q\alpha_{j+}$,
$c=\sum_{j=2}^q\alpha_{j-}$, $d=s_+$ and $e=s_{-}$. Condition $\rho(\Gamma)<1$ is equivalent to $\gamma<1$.
It is then clear that {\bf B1-B3} follow from {\bf G1-G2}. 

For the second part, it is easily seen that {\bf B4} is satisfied for a constant function $S$, $p=1$ and $s=\infty$. This gives the bound for $\Theta_{1,h}(H)$. If $\delta\geq 1$, the last bound for $\theta_{\delta,t}(Y)$ can be obtained from the inequalities
$$\vert x-y\vert^{\delta}\leq \left\vert x^{\delta}-y^{\delta}\right\vert,\quad x,y\geq 0.$$

\subsection{Proof of Proposition \ref{PARX}}

\begin{enumerate}
    \item To show the first point, we check the assumptions of Theorem \ref{autoreg}.
    We set $\zeta_t=\left(N^{(t)},Z_{t-1}\right)$, $E=\N\times \R_+$ and the state space $E$ is endowed with the $\ell_1-$norm.
		We first note that $(Y_t)_{t\in\Z}$ is a stationary solution of (\ref{Poisss}) if and only if
    $X_t = (Y_t, \lambda_t)'$ is the solution of $X_t=F\left(X_{t-1},\ldots,X_{t-q},\zeta_t\right)$
		with 
		$$F\left(x_1,\ldots,x_q,\zeta_t\right)=\left(N^{(t)}_{f(x_1, \ldots, x_q, Z_{t-1})},f(x_1, \ldots, x_q, Z_{t-1}) \right)'$$
	and $f(x_1, \ldots, x_q, Z_{t-1}) = \beta_0+\sum_{j=1}^q\beta_j s_j+\sum_{j=1}^q\alpha_j y_j+\pi'Z_{t-1}$, $x_i=(y_i,s_i)$, $1\leq i\leq q$.
 For $x \in (\N \times \R_+)^q$, 
      $$
      \E\left[\left\vert F\left(x,\zeta_1\right)\right\vert\right] = 2 \left(\beta_0+\sum_{j=1}^q\beta_j s_{j}+\sum_{j=1}^q\alpha_j y_{j}+\pi'\E(Z_1)\right) < \infty
      $$
      since $\E(\vert Z_1\vert) < \infty.$

 We then have, for $(x,x') \in {((\N \times \R_+)^q)}^2$ with $x = (x_1, \ldots, x_q), x' = (x_1', \ldots, x_q')~\forall j = 1, \ldots, q, x_j = (y_j, s_j), x_j = (y_j', s_j')$, 
 $$\E\left[\left\vert F(x,\zeta_t)-F(x',\zeta_t)\right\vert_{vec}\vert\mathcal{F}_{t-1}\right]\preceq \sum_{j=1}^q \begin{pmatrix}
\alpha_j & \beta_j \\
\alpha_j & \beta_j 
 \end{pmatrix} \left\vert x_j-x'_j\right\vert_{vec}.$$
In the previous bounds, we have used the identity $\E\left[\vert N^{(t)}_{h_{t-1}}-N^{(t)}_{g_{t-1}}\vert\big\vert \mathcal{F}_{t-1}\right]=\vert h_{t-1}-g_{t-1}\vert$ which is valid for two nonnegative random variables $h_{t-1},g_{t-1}$ measurable with respect to $\mathcal{F}_{t-1}$. The previous equality follows from the properties of the Poisson process.
Letting $$
\Gamma = \sum_{j=1}^q \begin{pmatrix}
\alpha_j & \beta_j \\
\alpha_j & \beta_j 
 \end{pmatrix},
$$
the matrix $\Gamma$ has two eigenvalues: $0$ and $\gamma$. Assumption {\bf PA1} then guarantees that $\rho(\Gamma)<1$. Assumptions {\bf B1-B3} of Theorem \ref{autoreg} are satisfied. Hence, according to Theorem \ref{autoreg}, there  exists a unique stationary and non-anticipative process $(X_t)_{t\in\Z}$ solution of \eqref{Poisss} and such that $\E\left[\vert X_t\vert\right]<\infty$. This process is ergodic. This completes the proof of the first point.

\item
For the second point, we use Proposition \ref{FDM}. To this end, it is only necessary to check {\bf B4} for $p=r=1$ and $s=\infty$. This is straightforward since we have the equality
$$\E\left[\vert F(x_1,\ldots,x_q,z,N^{(t)})-F(x_1,\ldots,x_q,\overline{z},N^{(t)})\vert\right]=2\vert \pi'(z-\overline{z})\vert.$$ 
The proof of the second point is now complete.

\end{enumerate}

\subsection{Proof of Proposition \ref{supercond}}
\begin{enumerate}
\item
We apply Theorem \ref{bibli}. To this end,
we define the random map from $E=\{0,1\}^q$ to $E$ by
$$f_t(x)=\left(\mathds{1}_{g(x,\zeta_t)>0},x_1,\ldots,x_{q-1}\right)'.$$
We set 
$$\delta_t:=\max_{y,y'\in \{0,1\}^q}\left\vert \mathds{1}_{g(y,\zeta_t)>0}-\mathds{1}_{g(y',\zeta_t)>0}\right\vert \leq \mathds{1}_{\max_{y\in\{0,1\}^q}g(y,\zeta_t)>0}-\mathds{1}_{\min_{y\in\{0,1\}^q} g(y,\zeta_t)>0}.$$
Setting $\left(y_t(x),\ldots,y_{t-q+1}(x)\right)'=f_1^t(x)$ for $t\geq q$, we have
$$y_t(x)=1\mbox{ if and only if } g\left(y_{t-1}(x),\ldots,y_{t-q}(x),\zeta_t\right)>0.$$ 
We have, setting $c_t=\left\vert y_t(x)-y_t(x')\right\vert$, 
$$c_t\leq \delta_t \max_{1\leq j\leq q}c_{t-j}.$$
Using the fact that $\delta_t\leq 1$, a straightforward induction on $i=0,\ldots,q-1$ shows that 
$$c_{t+i}\leq \delta_{t+i}\max_{1\leq j\leq q}c_{t-j}.$$
Setting $d(y,y')=\max_{1\leq i\leq q}\vert y_i-y_i'\vert$, 
$$d\left(f_1^{t+q-1}(x),f_1^{t+q-1}(x')\right)\leq \max_{0\leq i\leq q-1}\delta_{t+i} d\left(f_1^{t-1}(x),f_1^{t-1}(x')\right).$$
Setting $t=1$, this shows in particular that 
\begin{equation}\label{pre}
c\left(f_1^q\right)\leq \max_{1\leq i\leq q}\delta_i\leq \mathds{1}_{\max_{y,i}g(y,\zeta_i)>0} 
-\mathds{1}_{\min_{y,i}g(y,\zeta_i)>0}.
\end{equation}
From our assumptions, the last upper bound can vanish with positive probability,  
and then $\E\left[\log c(f_1^q)\right]=-\infty=\chi$.
Theorem \ref{bibli} leads to the result.

\item
The result will follow from Proposition \ref{WD}. To this end, we check Assumptions {\bf A1-A3}. {\bf A1} is automatic. We use the metric $d$ on $\{0,1\}^q$ which is bounded. We set $h_t=\max_{t-q+1\leq i\leq t}\delta_i$. Note that in the linear case, we have
$$h_t\leq \mathds{1}_{\phi_{+}+\max_{t-q+1\leq i\leq t}\upsilon_i>0}-\mathds{1}_{\phi_{-}+\min_{t-q+1\leq i\leq t}\upsilon_i>0}.$$
To check {\bf A2}, we use the bound (\ref{pre}) and the inequality $h_s\leq 1$ for all $s\in\Z$ to get
$$c\left(f_{t-Jq+1}^t\right)\leq \prod_{j=0}^{J-1} h_{t-jq}\leq h_t.$$
Moreover using (\ref{deJ}), we have
\begin{eqnarray*}
\E\left(h_t\vert \mathcal{F}_{t-Jq}\right)&\leq& 1-\P\left(\phi_{+}+\max_{t-q+1\leq i\leq t}\upsilon_i\leq 0\vert\mathcal{F}_{t-Jq}\right)-\P\left(\phi_{-}+\min_{t-q+1\leq i\leq t}\upsilon_i>0\vert\mathcal{F}_{t-Jq}\right)\\
&\leq& 1-\delta,
\end{eqnarray*}
provided that $Jq\geq q+K$. This guarantees {\bf A2}, with $\kappa=1-\delta$ and $m=Jq$.

Finally, let us check {\bf A3}. Setting $\overline{Z}_t=H\left(\eta_t,\ldots,\eta_1,\eta'_0,\eta_{-1},\ldots\right)$.
We have 
\begin{eqnarray*}
\E\left[d\left(f_t(y),\overline{f}_t(y)\right)\big\vert \mathcal{F}_{t-1}\vee\sigma(\xi_0')\right]&\leq&
\left\vert F_{\varepsilon}\left(-\sum_{i=1}^q a_i y_i-\pi'Z_{t-1}\right)-F_{\varepsilon}\left(-\sum_{i=1}^qa_i y_i-\pi'\overline{Z}_{t-1}\right)\right\vert\\
&\leq& L_{\varepsilon}\cdot\max_{j=1}^d\vert\pi_j\vert\cdot \sum_{j=1}^d\left\vert Z_{j,t-1}-\overline{Z}_{j,t-1}\right\vert,
\end{eqnarray*}
where $L_{\varepsilon}$ denotes the Lipschitz constant of $F_{\varepsilon}$.
\end{enumerate}

\subsection{Proof of Proposition \ref{cat}}
\begin{enumerate}
\item
For the first point, we will apply Theorem \ref{bibli} for the discrete metric $d(x,y)=\mathds{1}_{x\neq y}$.
Setting for $x\in E^q$ and $t\in\Z$, $f_{t-q+1}^t(x)=f_t\circ\cdots\circ f_{t-q+1}(x)$, we have for $x,y\in E^q$, 
\begin{eqnarray*}
\left\{f_{t-q+1}^t(x)=f_{t-q+1}^t(y)\right\}&\supset&\left\{f_{t-q+1}^t(x)=f_{t-q+1}^t(y)=(1,\ldots,1)\right\}\\
&=&\left\{\varepsilon_j\in [0,\eta(Z_{j-1})]: j=t-q+1,\ldots t\right\}.
\end{eqnarray*}
We then obtain
$$d\left(f_{t-q+1}^t(x),f_{t-q+1}^t(y)\right)\leq \left(1-\prod_{j=t-q+1}^t\mathds{1}_{\varepsilon_j\in [0,\eta(Z_{j-1})]}\right)d(x,y).$$
Let us show that 
\begin{equation}\label{amelior}
p:=\P\left(\varepsilon_j\in [0,\eta(Z_{j-1})]: j=t-q+1,\ldots t\right)>0.
\end{equation}
Note that by stationarity, $p$ does not depend on $t$. Assume that $p=0$. 
From Assumption {\bf C2} and the properties of the conditional expectations, we have 
$$p=\E\left[\eta(Z_{t-1})\prod_{j=t-q+1}^{t-1}\mathds{1}_{\left\{\varepsilon_j\leq \eta(Z_{j-1})\right\}}\right]=0.$$
Since $\eta$ is positive, we get $\P\left(\varepsilon_j\in [0,\eta(Z_{j-1})]: j=t-q+1,\ldots t-1\right)=0$.
By finite induction, we deduce that $\P\left(\varepsilon_{t-q+1}\leq \eta(Z_{t-q})\right)=0$. Since this latter probability equals to $\E\left(\eta(Z_{t-q})\right)>0$, we obtain a contradiction. The property (\ref{amelior}) is then valid and Theorem
\ref{bibli} applies, which leads to the conclusion.
\item
For the second point, we will use Theorem \ref{main} and Proposition \ref{WD}. 
Since $d$ is a bounded metric, {\bf A1} is automatically satisfied.
Next, observe that 
\begin{eqnarray*}
\P\left(f_{t-q+1}^t(x)=f_{t-q+1}^t(y)\vert \mathcal{F}_{t-q}\right)&\geq& \P\left(f_{t-q+1}^t(x)=f_{t-q+1}^t(y)=(1,\ldots,1)\vert\mathcal{F}_{t-q}\right)\\
&\geq& \P\left(\varepsilon_t,\ldots,\varepsilon_{t-q+1}\in [0,\eta_{-}]\mathcal{F}_{t-q}\right)\\
&=&\P\left(\varepsilon_t,\ldots,\varepsilon_{t-q+1}\in [0,\eta_{-}]\right)\\
&\geq& \eta_{-}^q.
\end{eqnarray*}
This yields to the bound $\E\left[d\left(f_{t-q+1}^t(x),f_{t-q+1}^t(y)\vert \mathcal{F}_{t-q}\right)\right]\leq 1-\eta_{-}^q$, which shows the second part of {\bf A2}. The first part is automatic. 

It remains to check {\bf A3'}. Note that for $i,j\in E$, we have $(N-1)^{-1}\vert i-j\vert\leq\mathds{1}_{i\neq j}\leq \vert i-j\vert$. Using the $\ell_1-$metric on $E^q$ which is equivalent to the discrete metric, we have 
\begin{eqnarray*}
\E\left[\vert K^{-}_z(\varepsilon_1\vert y_1,\ldots,y_q)-K^{-}_{z'}(\varepsilon_1\vert y_1,\ldots,y_q)\vert\right]&\leq&
\int_0^1\left\vert K^{-}_z(u\vert y_1,\ldots,y_q)-K^{-}_{z'}(u\vert y_1,\ldots,y_q)\right\vert du\\
&\leq&\sum_{j=1}^N\left\vert \sum_{i=1}^j K_z(i\vert y_1,\ldots,y_q)-\sum_{i=1}^j K_{z'}(i\vert y_1,\ldots,y_q)\right\vert\\
&\leq & NC\left\vert z-z'\right\vert.
\end{eqnarray*}
Assumption {\bf A3'} then follows with $s=\infty$ and $r=p=1$.
The bound for the functional dependence coefficients is then a direct consequence of Proposition \ref{WD}.$\square$
\end{enumerate}

\bibliographystyle{plainnat}
\bibliography{bibbib}

\end{document}